\documentclass[final]{siamltex1213}

\newcommand{\doingrevtex}[1]{}
\newcommand{\doingsiam}[1]{#1}

\usepackage{amsmath}
\usepackage{graphicx}
\usepackage{bmpsize}

\usepackage[usenames,dvipsnames]{xcolor}

\doingrevtex{\usepackage{amsthm}}

\doingsiam{\usepackage{breakcites}
\usepackage{microtype}
}

\usepackage{epstopdf}

\newcommand{\mat}[1]{\boldsymbol{#1}}

\newcommand{\ot}{  {\scriptstyle \otimes}_{ \tau } }
\newcommand{\ots}{ {\scriptstyle \otimes}_{ \! \tau_s } }
\newcommand{\oto}{ {\scriptstyle \otimes}_{ \! \tau_0 } }
\newcommand{\otone}{ {\scriptstyle \otimes}_{ \! \tau_1 } }
\newcommand{\otm}{ {\scriptstyle \otimes}_{ \! \tau_m } }

\newcommand{\otpm}{ {\scriptstyle \otimes}_{ \! \tau_{m+1}}}

\newcommand{\xpose}{ {\scriptscriptstyle T}}

\doingrevtex{
\newtheorem{thm}{\protect\theoremname}
\theoremstyle{plain}

\theoremstyle{remark}

\theoremstyle{plain}
\newtheorem{prop}[thm]{\protect\propositionname}
\providecommand{\lemmaname}{Lemma}
\providecommand{\propositionname}{Proposition}
\providecommand{\remarkname}{Remark}
\providecommand{\theoremname}{Theorem}
}

\doingsiam{
\title{A $N$-Body Solver for Square Root Iteration \thanks{This article was released under LA-UR-15-26304.  The Los Alamos National
Laboratory is operated by Los Alamos National Security, LLC for the NNSA of the
USDoE under Contract No.  DE-AC52- 06NA25396.}}
\author{Matt Challacombe\footnotemark[1] \footnotemark[2],
        Terry Haut\footnotemark[1]~and
        Nicolas Bock\footnotemark[1]  }
}

\begin{document}
\doingrevtex{
\title{A $N$-Body Solver for Square Root Iteration}
\author{Matt Challacombe, Terry Haut and Nicolas Bock}
}

\doingsiam{
\maketitle
\renewcommand{\thefootnote}{\fnsymbol{footnote}}
\footnotetext[1]{Los Alamos National Laboratory}
\footnotetext[2]{matt.challacombe@freeon.org}
\renewcommand{\thefootnote}{\arabic{footnote}}

\pagestyle{myheadings}
\thispagestyle{plain}
\markboth{Challacombe, Haut \& Bock}{A $N$-Body Solver for Square Root Iteration}
}

\doingrevtex{
\affiliation{Los Alamos National Laboratory}
\preprint{\tt LA-UR-15-26304}
}

\begin{abstract}
We develop the Sparse Approximate Matrix Multiply ($\tt SpAMM$) $n$-body solver for first order Newton Schulz iteration of the
matrix square root and inverse square root.
The solver performs recursive two-sided metric queries on a modified Cauchy-Schwarz criterion,
culling negligible sub-volumes of the product-tensor for problems with structured decay in the sub-space metric.
These sub-structures are shown to bound the relative error in the matrix-matrix product, and
in favorable cases, to enjoy a reduced computational complexity governed by dimensionality reduction of
the product volume.  A main contribution is demonstration of a new, algebraic locality that develops under
contractive identity iteration, with collapse of the metric-subspace onto the identity's plane diagonal,
resulting in a stronger $\tt SpAMM$ bound.  Also, we carry out a first order {Fr\'{e}chet} analyses for single and
dual channel instances of the square root iteration, and look at bifurcations due to ill-conditioning and a
too aggressive $\tt SpAMM$ approximation.  Then, we show that extreme $\tt SpAMM$ approximation and contractive
identity iteration can be achieved for ill-conditioned systems through regularization, and we
demonstrate the potential for acceleration with a scoping, product representation of the inverse factor.
\end{abstract}

\maketitle
\section{Introduction}
In many areas of current numerical interest, matrix equations with decay properties describe correlations over a range of scales.
By decay, we mean an approximate inverse relationship between a matrix element's magnitude and an associated distance;
this might be a slow inverse exponential relationship between matrix elements and a Cartesian separation,
or it might involve a non-Euclidean distance, {\em  e.g.}~between character strings.

A common approach to exploiting matrix decay involves sparse approximation of inverse factors that transform Gramian
equations to a representation independent form, via congruence transformations based on {L\"{o}wdin's} symmetric orthogonalization
(the matrix inverse square root) \cite{Lowdin56,Naidu11}, inverse Cholesky factorization \cite{krishtal2015} or related transformations that involve an inverse
or pseudo-inverse \cite{head2000tensors,bjorck2004calculation,grohs2013intrinsic,grohs2014intrinsic}.
Gramian inverse factors with decay are ubiquitous to problems with local, non-orthogonal support, including
finite element calculations \cite{demko1984decay,Haber2014}, radial-basis-function finite-difference calculations \cite{tolstykh2003using,shankar2016radial},
in the ``direct'' approach to radial-basis interpolation \cite{sarra2014}, with frames \cite{fornasier2005intrinsic,heil2007history},
with computation involving  ``lets'' of various types \cite{grohs2013intrinsic,grohs2014intrinsic},
and in the Atomic Orbital (AO) representation \cite{Jansik2007,helgaker2008molecular}.

\begin{figure}[h]
\doingrevtex{\includegraphics[width=3.5in]{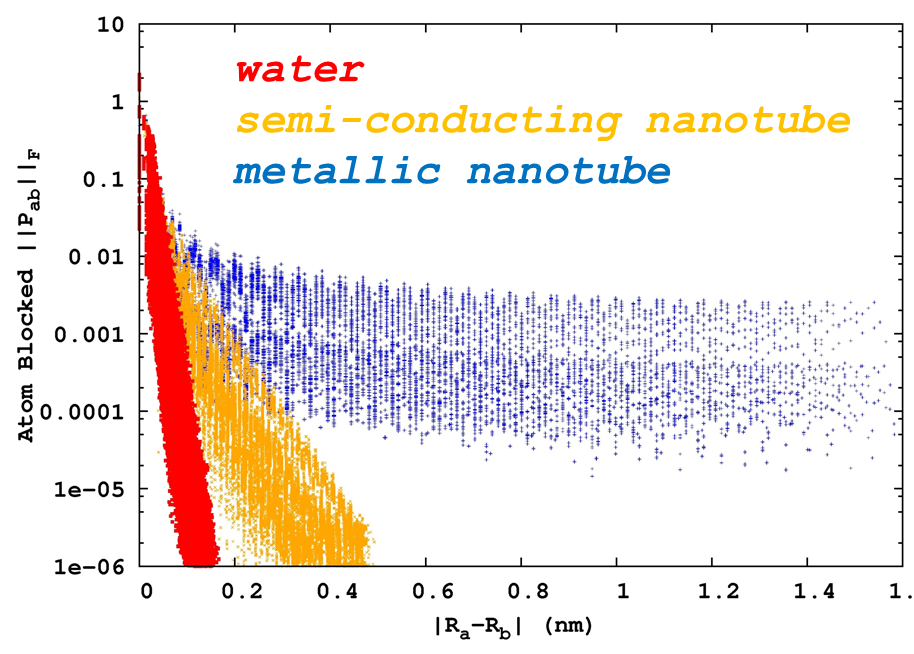}}
\doingsiam{\includegraphics[width=5.0in]{decay_picture.png}}
  \caption{\label{figure1} Examples from electronic structure of decay for the
    spectral projector (gap shifted sign function) with respect to
    the local (atomic) support.  Shown is decay for systems with
    correlations that are short (insulating water), medium
    (semi-conducting 4,3 nanotube), and long (metallic 3,3 nanotube)
    ranged, from exponential (insulating) to algebraic (metallic). }
\end{figure}

Off-diagonal decay of the matrix sign function is also a well developed area of study in statistics and statistical physics
\cite{penrose1974,voit00, Anselin2003, Hardin2013, Krishtal2014}, and in electronic structure, where sparse approximation enables
fast computation of the the gap shifted matrix sign function as projector of an effective Hamiltonian
\cite{Benzi99,Benzi02,Bowler12,Benzi13}.  Short to long ranged decay properties of the projector are shown in Fig.~\ref{figure1}.
These matrix functions, the matrix sign function and the matrix inverse square root, are related by Higham's identity \cite{Higham08}:
\begin{equation}\label{highamsid}
\rm{sign} \left( \begin{bmatrix} 0 & \mat{s}      \\ \mat{I}       & 0\end{bmatrix} \right)  =
                 \begin{bmatrix} 0 & \mat{s}^{1/2} \\ \mat{s}^{-1/2} & 0\end{bmatrix}  .
\end{equation}

A well conditioned matrix $\mat{s}$ may often correspond to matrix
sign and inverse square root functions with rapid exponential decay,
and be amenable to {\em ad hoc} matrix truncation or ``sparsification'',
$\bar{\mat{s}} = \mat{s}+ \mat{\epsilon}^{\mat{s}}_\tau$,
where $\mat{\epsilon}^{\mat{s}}_\tau$ is the error introduced
according to some criterion $\tau$, supported by useful bounds
to matrix function elements \cite{Benzi1999,Benzi:2007:Decay,Passenbrunner14,Haber2014,canuto2014decay}.
The criterion $\tau$ might be a drop-tolerance,
$\epsilon^{\mat{s}}_{\tau} = \{-s_{ij}*\hat{\mat{e}}_i \, | \, |s_{ij}|<\tau \}$,
a radial cutoff,
$\epsilon^{\mat{s}}_{\tau} = \{-s_{ij}*\hat{\mat{e}}_i \, | \, \lVert \mat{r}_i - \mat{r}_j \rVert > \tau \}$,
or some other approach to truncation, perhaps involving a sparsity
pattern chosen {\em a priori} for computational expedience.
Then, the sparse general matrix-matrix multiply ($\tt{SpGEMM}$)
\cite{Gustavson78, Toledo97, challacombe00, bowler00} may be employed, yielding fast
solutions for multiplication rich iterations, and  with fill-in modulated by truncation.
Exhaustive surveys of these methods in the numerical linear algebra are given by Benzi
\cite{Benzi99,Benzi02}, and by Bowler \cite{Bowler12} and Benzi \cite{Benzi13}
for electronic structure.

In addition to sparsity, data localities leading to high operation counts are essential
for kernels like the $\tt SpGEMM$ and their distributed implementations.
Over the past decades, methods have evolved from bandwidth reduction (Cuthill-McKee) + greedy blocking \cite{Toledo97}, progressing
with tours of the graph via heuristic solutions to the Traveling Salesman Problem (TSP) \cite{Pinar1999,Akbudak2012,Lemire2012},
and more recently towards {reordering} based on cache modeling and dynamic sampling \cite{Frasca2012,Pichel2014}.
{Ordering} with graph {partitioning}, targeting the load balance, may also lead to {exploitable} localities,
via {\em e.g.}~proximity to the diagonal \cite{buluc2013recent}.  Of current interest are ordering schemes that enhance the
{weighted} block-locality of the Page Rank problem \cite{kamvar2003,DelCorso2005,li2009,yan2014}.

\begin{figure}[h]
  \begin{center}
    \fbox{\includegraphics[width=0.47\columnwidth]{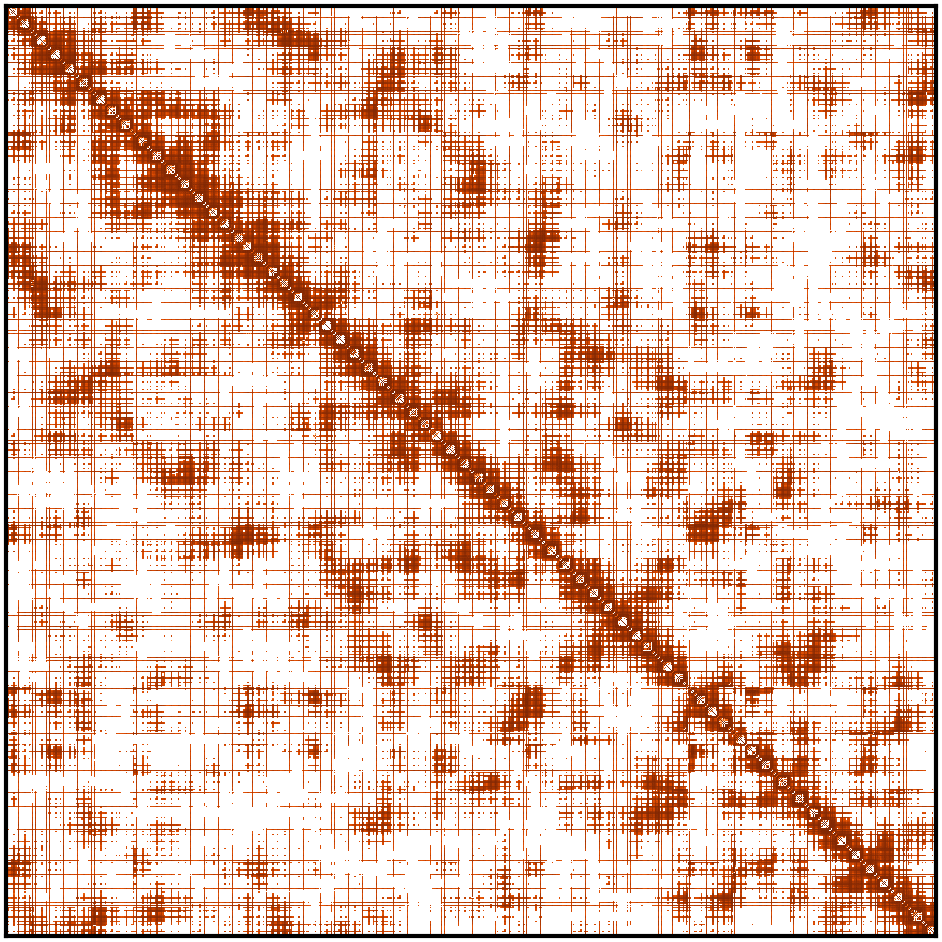}}
    \fbox{\includegraphics[width=0.47\columnwidth]{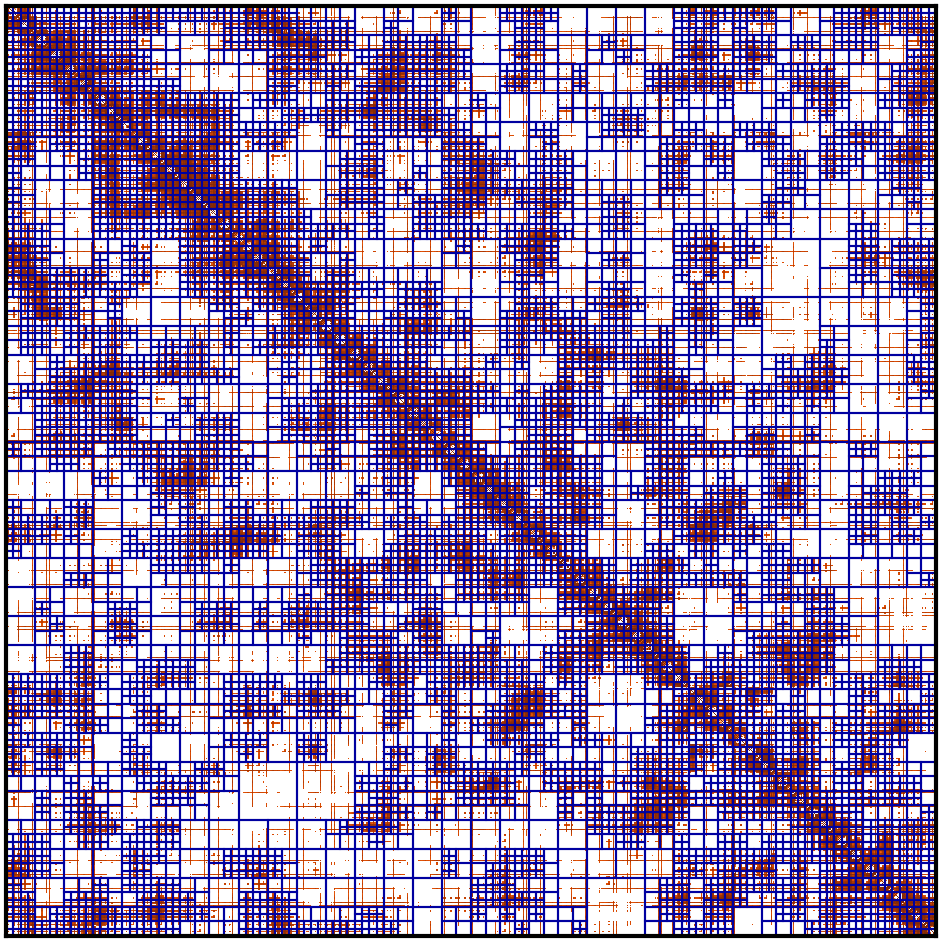}}
  \end{center}
  \caption{\label{nickspick}
    At left, Block-By-Magnitude (BBM) structure of a quantum chemical Gramian (overlap matrix),
    for a box of 100 water molecules, with Cartesian support along a locality preserving curve.
    At right, quadtree resolution of neighborhoods with norms down to $10^{-3}$.}
\end{figure}

Matrix locality may also result from {an ordering that preserves} locality in an {auxiliary} representation,
a property of sub-space mappings that preserve local neighborhoods \cite{Belkin2002,Belkin2003,Belkin2008}.
In the case of electronic structure,  Space Filling Curve (SFC) heuristics applied
to a local Cartesian basis results in Gramian matrices with neighborhoods segregated by magnitude \cite{challacombe00,brazdova08}, as shown in Fig.~(\ref{nickspick}).
Likewise, Sierpinski curves and Self Avoiding Walks on meshes lead to locality preserving orderings \cite{Heber1998,bader13},
for {\em e.g.}~finite elements \cite{Oliker2002,Schamberger04}.
This type of weighted block-locality or ``Block-By-Magnitude'' (BBM) structure of the subspace metric $\lVert \cdot \rVert_F$
is finely resolved with the quadtree matrix
\cite{Wise:1984:RMQ:1089389.1089398,springerlink:10.1007/3-540-51084-2_9,
      Samet:1990:DAS:77589,Wise1990,Wise:Ahnentafel,
      Lorton:2006:ABL:1166133.1166134};
\begin{equation}
\mat{a}^i = \begin{bmatrix} \,  \mat{a}^{i+1}_{00} \, & \,  \mat{a}^{i+1}_{01} \,
\\[0.2cm]  \, \mat{a}^{i+1}_{10} \,  & \,\mat{a}^{i+1}_{11} \, \end{bmatrix} \, ,
\end{equation}
where $i$ is the recursion depth, and
\begin{equation}
\lVert \mat{a}^i \rVert_F = \sqrt{  \lVert \mat{a}^{i+1}_{00}\rVert^2_F + \lVert \mat{a}^{i+1}_{01}\rVert^2_F
                                  +  \lVert \mat{a}^{i+1}_{10}\rVert^2_F + \lVert \mat{a}^{i+1}_{11}\rVert^2_F } \; ,
\end{equation}
is the sub-multiplicative Frobenius norm \cite{Golub96,Higham02,Kahan2013}.

However, despite structuring for cache, distributed memory or to enhance BBM structuring, matrix truncation may still be ineffective for
ill-conditioned problems, because the rate of decay may be too slow, and also because of increased numerical sensitivities
to the sparse approximation:
\begin{equation} \label{sparseapprox}
\overline{ \mat{a} \cdot \mat{b} }\; = \; \mat{a}\cdot\mat{b} \; +\; \mat{\epsilon}^{\mat{a}}_\tau \cdot \mat{b} \;+\;
 \mat{a} \cdot \mat{\epsilon}^{\mat{b}}_\tau  \; + \;   {\mathcal O}(\tau^2) \, ,
\end{equation}
allowing to control only absolute errors.
An alternative approach is to find a reduced rank approximation, ideally closed under the operations of interest.
However, rank reduction may be expensive if the rank is not much, much smaller than the dimension.
Interestingly, in the ultra-flat limit, kernel methods enjoy rank reduction corresponding formally to change of
basis, enabling fast methods for constructing the generalized inverse \cite{schaback2007convergence,Challacombe2016}.
In cases with simply slow exponential decay however, our experience has so far been that naive element dropping is
about as effective as dropping singular values.

In this contribution, we consider the regime between trivial sparsity and formal rank reduction, with
fast multiplications exploiting instead an accelerated {\em volumetric decay} in subspace-metric of the product-tensor.
There, the $\tt SpAMM$ kernel {carries} out octree scoping of low-dimensional structures that bound
the relative error while yielding a reduced complexity multiplication.   Beyond decay associated with the
the matrix square root and its inverse, we demonstrate additional compression of these bounding volumes under
contractive identity iteration.

This paper is organized as follows:  In Section \ref{spammsection},  we modify the $\tt SpAMM$ occlusion-cull to
bound the relative product error.  In Section \ref{nsiterations}, we review
several instances of the first order {Newton-Schulz (NS)} square root iteration, and go over the contractive
identity iteration that develops in the basin of stability.  In Section \ref{implementation}, we overview
a generic implementation of the $\tt SpAMM$ kernel and introduce quantum chemical and {engineering} data of
interest.  In Section \ref{errorflow}, we develop a {Fr\'{e}chet} {analysis} for NS instances and the $\tt SpAMM$ algebra,
and examine error flows in bifurcating and stable square root iterations for ill-conditioned problems.
In Section \ref{regularization} we show that even difficult, ill-conditioned problems can be brought to the
regime of strongly contractive identity iteration, through iterative regularization and precision scoping.
In Section \ref{locality}, we show for the first time the process of lensing, involving sub-space contraction to
diagonal planes of the $ijk$-cube ($i=j$, $i=k$ and/or $j=k$), followed finally by compression onto the identity's plane diagonal,
{yielding} additional orders of magnitude compression of $\tt SpAMM$ sub-volumes.  Finally, in Section \ref{conclusions}
we argue it may be possible to remain close to the lensed state whilst constructing
a {deferred product representation} of the inverse factor.


\section{$\tt SpAMM$}\label{spammsection}
The Sparse Approximate Matrix Multiply ($\tt SpAMM$) is a reduced complexity approximation that
evolved from a row-column skipout mechanism within the blocked-compressed-sparse-row (BCSR) \cite{challacombe1999simplified} and
the distributed-blocked-compressed-row (DBCSR) data structures \cite{Challacombe:2000:SpMM},
to methods with fast subspace resolution through octree {recursion} \cite{Challacombe2010,Bock2013,BockCK14}.
Finding sub-spaces via fast range or metric {query} is a generic $n$-body problem handled with agility by the
quadtree \cite{genericityindata,Geerts2002,Samet:2006:DBDS,Gottschling2009}, a problem related to
spatial hashing \cite{Teschner2003,Lefebvre2006} and the occlusion-cull in visualization \cite{Pantazopoulos2002}.

The $\tt SpAMM$ kernel $\ot$ provides fast approximate multiplication for matrices with decay and metric locality,
with errors {controlled} by the scoping parameter $\tau$:
\begin{equation}\label{approxproduct}
\widetilde{\mat{a}\cdot \mat{b}} \,  \equiv \, \mat{a} \ot \mat{b} \,
  = \, \mat{a} \cdot \mat{b} + \mat{\Delta}^{a \cdot b}_{\tau} \, .
\end{equation}
As $\tau \rightarrow 0$, $\tt SpAMM$ reverts to the recursive $\tt GEMM$ \cite{Gustavson99,Elmroth2004}.

In this work, we promote the following stable version of the $\tt SpAMM$ occlusion-cull:
\doingrevtex{\begin{widetext}}
\begin{equation}\label{newspamm}
\mat{a}^{i} \ot \mat{b}^{i} =
\left\{
        \begin{array}{ll}
                 \mat{0} \quad \tt{if}\quad \lVert \mat{a}^i \rVert \lVert \mat{b}^i \rVert < \tau \lVert a \rVert \lVert b \rVert \\[0.2cm]
                 \mat{a} ^i \cdot \mat{b}^i \quad  \tt{if}(i=\tt{leaf}) \\[0.2cm]
\begin{bmatrix} \mat{a}^{i+1}_{00} \ot \mat{b}^{i+1}_{00} +\mat{a}^{i+1}_{01} \ot \mat{b}^{i+1}_{10} \; , \; &
                \mat{a}^{i+1}_{00} \ot \mat{b}^{i+1}_{01} +\mat{a}^{i+1}_{01} \ot \mat{b}^{i+1}_{11}  \\[0.2cm]
                \mat{a}^{i+1}_{00} \ot \mat{b}^{i+1}_{01} +\mat{a}^{i+1}_{01} \ot \mat{b}^{i+1}_{11} \; , \; &
                \mat{a}^{i+1}_{00} \ot \mat{b}^{i+1}_{01} +\mat{a}^{i+1}_{01} \ot \mat{b}^{i+1}_{11}
\end{bmatrix}  \quad \tt{else}
                \end{array}
              \right.  \, ,
\end{equation}
\doingrevtex{\end{widetext}}
with $\lVert \cdot \rVert \equiv \lVert \cdot \rVert_F$ and the leaf condition determined by the block size, $N_b$.
This scoping partitions the product tensor into two sub-spaces:  the space of culled leaf-tasks, $\mat{a} \ot \mat{b}$,
and its complement, the occlusion error $\mat{\Delta}^{a \cdot b}_{\tau}$ of avoided multiplications.
This occlusion error is bounded by
\begin{equation}\label{bound}
\frac{\lVert \mat{\Delta}^{a \cdot b}_{\tau} \rVert}{\lVert \mat{a} \rVert  \,  \lVert \mat{b} \rVert }  \, \leq \, n^2 \tau \, ,
\end{equation}
as shown in the following section, a result commensurate with the stable, normwise multiplication criteria emphasized
by Demel, Dumitriu, Holtz and Kleinberg (DDHK) \cite{Demmel:2007:FastMM}.

\subsection{Bound}
We now prove Eq.~(\ref{bound}):\\
\doingsiam{\begin{proposition}\label{lem:SpAMM mult, prop}}
\doingrevtex{\begin{prop}}
Let $\tau_{ab} = \tau \| \mat{a} \| \| \mat{b} \| $. Then for each $i,j$,
\[
\left|\left(\mat{a}\ot \mat{b}\right)_{ij}-\left( \mat{a} \cdot \mat{b} \right)_{ij}\right| \, \leq \, n \,  \tau_{ab},
\]
and
\[
\left\Vert \mat{a} \ot \mat{b}- \mat{a} \cdot \mat{b} \right\Vert \, \leq  \, n^{2} \,\tau_{ab}.
\]
\doingsiam{\end{proposition}}
\doingrevtex{\end{prop}}

\begin{proof}
We first show the following technical result: it is possible to choose $\alpha_{lij}\in\left\{ 0,1\right\} $ such that
\begin{equation}
\left( \mat{a} \ot \mat{b}\right)_{ij}=\sum_{l=1}^{n}a_{il}\, b_{lj} \, \alpha_{lij},\label{eq:spamm form, lemma}
\end{equation}
In addition, if $\alpha_{lij}=0$, then \textup{$\left|a_{il}\right|\left|b_{lj}\right|<\tau_{ab}$}. To show this, we use
induction on the number $k_{\max}$ of levels.

First, if $k_{\max}=0$,
\[
\mat{a} \ot \mat{b}=\begin{cases}
0 & \,\,\text{if}\,\,\left\Vert \mat{a} \right\Vert \left\Vert \mat{b} \right\Vert <\tau_{ab},\\
\mat{a} \cdot \mat{b} & \,\,\text{else}.
\end{cases}
\]
 Therefore, $\mat{a}\ot \mat{b}$ is of the form (\ref{eq:spamm form, lemma})
with either all $\alpha_{lij}=0$ or all $\alpha_{lij}=1$. Moreover,
if $\alpha_{lij}=0$, then $\left|a_{il}\right|\left| b_{lj}\right|\leq\left\Vert \mat{a} \right\Vert
 \left\Vert \mat{b}\right\Vert <\tau_{ab}$.

Now assume that the claim holds for $k_{\max}-1$. We show that it
holds for $k_{\max}$. Indeed, if $\left\Vert \mat{a} \right\Vert \left\Vert \mat{b} \right\Vert < \tau_{ab}$,
we have that $\mat{a} \ot \mat{b}=0$, which is of the form (\ref{eq:spamm form, lemma})
with all $\alpha_{lij}=0$. Also, if $\alpha_{lij}=0$, then $\left| a_{il}\right|\left|b_{lj}\right|
<\left\Vert \mat{a}\right\Vert \left\Vert\mat{b}\right\Vert <\tau_{ab}$.

Now assume that $\left\Vert \mat{a}\right\Vert \left\Vert \mat{b}\right\Vert \geq\tau_{ab}$.
Then
\[
\mat{a} \ot \mat{b}=\left(\begin{array}{cc}
\mat{a}_{00}\ot \mat{b}_{00}+\mat{a}_{01}\ot \mat{b}_{10} & \mat{a}_{00}\ot \mat{b}_{01}+\mat{a}_{01}\ot \mat{b}_{11}\\
\mat{a}_{10}\ot \mat{b}_{00}+\mat{a}_{11}\ot \mat{b}_{10} & \mat{a}_{10}\ot \mat{b}_{10}+\mat{a}_{11}\ot \mat{b}_{11}
\end{array}\right).
\]
We need to consider four cases: $i\leq n/2$ and $j\leq n/2$, $i>n/2$
and $j>n/2$, $i>n/2$ and $j\leq n/2$, and, finally, $i>n/2$ and
$j>n/2$. Since the analysis is similar for all four cases, we only
consider $i\leq n/2$ and $j\leq n/2$. We have that
\begin{eqnarray*}
\left( \mat{a} \ot \mat{b}\right)_{ij} & = & \left( \mat{a}_{00} \ot \mat{b}_{00}+\mat{a}_{01}\ot \mat{b}_{10}\right)_{ij}\\
 & = & \sum_{l=1}^{n/2}\left(\mat{a}_{00}\right)_{il}\left(\mat{b}_{00}\right)_{lj}\alpha_{lij}^{0}  \\
 &   & \qquad + \sum_{l=1}^{n/2}\left(\mat{a}_{01}\right)_{il}\left(\mat{b}_{10}\right)_{lj}\alpha_{lij}^{1}\\
 & = & \sum_{l=1}^{n}a_{il}b_{lj}\alpha_{lij},
\end{eqnarray*}
where we used the induction hypothesis in the second equality.

Now suppose that $\alpha_{lij}=0$ for some $l$. Then $\tilde{\alpha}_{lij}^{0}=0$
if $l\leq n/2$ or $\tilde{\alpha}_{l-n/2,ij}^{1}=0$ if $l>n/2$.
If, e.g., $\tilde{\alpha}_{l-n/2,ij}^{1}=0$, then
$\left|a_{il}\right|\left|b_{lj}\right|=\left|\left(\mat{a}_{01}\right)_{i,l-n/2}\right|\left|\left(\mat{b}_{10}\right)_{l-n/2,j}\right|<\tau_{ab}$,
where we used the induction hypothesis in the final inequality. The
analysis for $l\leq n/2$ is similar, and the claim
follows.

We can now finish the proof of Proposition~\ref{lem:SpAMM mult, prop}. Indeed, by (\ref{eq:spamm form, lemma}),
\begin{eqnarray*}
\left|\left( \mat{a} \ot \mat{b}\right)_{ij}-\left( \mat{a} \cdot \mat{b} \right)_{ij}\right|
& \leq & \sum_{l=1}^{n} \left| a_{il}b_{lj} \right|   \left| \alpha_{lij}-1 \right|   \\
& = & \sum_{\alpha_{lij}=0}\left|a_{il}b_{lj}\right|.
\end{eqnarray*}
In addition, if $\alpha_{lij}=0$, then $\left|a_{il}b_{lj}\right|<\tau_{ab}$
and the lemma follows.

\end{proof}

\subsection{Related research}\label{relatedr}

$\tt SpAMM$ is perhaps most closely related to the Strassen-like branch of fast matrix multiplication
\cite{springerlink:10.1007/BF02165411,DDH07,Ballard2014,LeGall:2014:PTF:2608628.2608664,Ambainis15},
and also methods for group theoretical embedding allowing fast polynomial multiplication
\cite{cohn2003group,cohn2005group,Umans:2006:GAM:1145768.1145772}.
In the Strassen-like approach, disjoint volumes in high order tensor expansions of the product are recursively excluded,
while in the $\tt SpAMM$ approach to fast multiplication, the subspace metric of the product tensor is
recursively queried for occlusion of negligible volumes, with error bounded by Eq.~(\ref{bound}).
These  methods for fast matrix multiplication are stable, satisfying the DDHK normwise product bound \cite{DDH07}.

This work offers a data local alternative to fast non-deterministic methods for sampling the product,
which include sketching \cite{Sarlos2006,Drineas2006,Mahoney2012,Pagh2013,Woodruff2015},
joining \cite{Mishra92,Hoel94,Jacox03,Chen07,Amossen09,Lieberman08,Kim09},
sensing \cite{iwen2009note} and probing \cite{chiu2012matrix}.  These  methods may involve probabilistic
assumptions and on the fly sampling, with the potential for complexity reduction due to statistical approximations.
$\tt SpAMM$  also employs on the fly weighted sampling,  with compression through octree scoping of
metric tensor decay, and with additional subspace compression due to the onset of identity iteration.

$\tt SpAMM$ is related to the generalized $n$-body methods popularized by {Gray} \cite{Gray01,Gray2003}.
Here and in related research, we are interested in generic approaches to approximation that are data agnostic,  based on the quadtree and its
generalizations \cite{genericityindata,Geerts2002,Samet:2006:DBDS,Gottschling2009}
and and on the facile measure $\lVert \cdot \rVert \equiv \lVert \cdot \rVert_F$ \cite{Kahan2013}.
In this work, a fast two-sided metric {query} enables octree scoping with the occlusion criteria
$\lVert \mat{a}^i \rVert \lVert \mat{b}^i \rVert < \tau \lVert a \rVert \lVert b \rVert$.
With quantum chemical Fock exchange, a fast three-sided metric {query} enables hextree scoping with
a related, Cauchy-Schwarz like occlusion criterion (direct SCF) \cite{Challacombe2014}.
It may also be possible to exploit subspace locality more broadly, though mappings that optimally preserve
local neighborhoods in higher dimensions, {\em e.g.} via the Laplace-Beltrami {operator} \cite{Belkin2002,Belkin2003,Belkin2008}.

For distributed architectures, $n$-body methods offer well established protocols for turning spatial locality into
data and temporal locality \cite{Warren:1993:PHO:169627.169640,Warren:1995:HOT,Warren1995266,WarrenGordonBell1997,Warren2013}.
Recently, we showed strong scaling for the $\tt SpAMM$ kernel \cite{BockCK14}, while
Driscoll {\em et.~al} were able to show perfect strong scaling and communication optimally
for pairwise $n$-body methods \cite{Driscoll2013}.  A uniform approach to generic scoping is empowered at the ecosystems
level by runtime support for recursive task parallelism
\cite{dinan2008scioto,sampath2008dendro,Lashuk:2009:SFC,meng2010dynamic,sampath2010parallel,min2011hierarchical,meng2014scalable}.

Finally, this work is inspired broadly by Higham's work, particularly by Higham, Mackey, Mackey and Tisseur (HMMT) in
2005 \cite{higham2005} on square root iteration and the group structure of matrix functions.
Also, it is influenced by Chen and Chow's \cite{chen2014} approach to scaled NS iteration for ill-conditioned problems, and by
the Helgaker group's work on NS iteration, whose notation we follow in part \cite{Jansik2007}.


\section{Newton-Shulz Iterations} \label{nsiterations}

There are two common, first order NS iterations; the sign iteration
and the square root iteration, related by the square $\mat{I}\left(
\cdot \right)= {\rm sign}^2\left( \cdot \right) $.  These equivalent
iterations converge linearly at first, then enter a basin of stability
marked by super-linear convergence.

\subsection{Sign iteration}

For the NS sign iteration, this basin is marked by a behavioral change
in the difference $\delta \mat{X}_k = \widetilde{\mat{X}}_k -\mat{X}_k
= {\rm sign} \left(\mat{X}_{k-1}+\delta \mat{X}_{k-1} \right) -{\rm
  sign} \left(\mat{X}_{k-1} \right)$, where $\delta \mat{X}_{k-1}$ is
some previous error.  The change in behavior is associated with the
onset of idempotence and the bounded eigenvalues of ${\rm sign}'\left(
\cdot \right)$, leading to stable iteration when ${\rm sign}'\left(
\mat{X}_{k-1} \right) \delta \mat{X}_{k-1} < 1 $.  Global perturbative
bounds on this iteration have been derived by Bai and Demmel
\cite{Bai98usingthe}, while Byers, He and Mehrmann \cite{byers1997matrix} developed
asymptotic bounds.  The automatic stability of sign iteration is a
well developed theme in Higham's work~\cite{Higham08}.

\subsection{Square root iteration}
We are concerned with resolution of the identity
\begin{equation}
\mat{I} \left( \mat{s} \right) =\mat{s}^{1/2} \cdot \mat{s}^{-1/2} \, ,
\end{equation}
and its low-complexity computation with fast methods.

Starting with eigenvalues rescaled to the domain $(0,1]$ with the easily obtained
largest eigenvalue,   $\mat{s} \leftarrow \mat{s}/s_{N-1}$, and with $\mat{z}_0=\mat{I}$ and
$\mat{x}_0=\mat{y}_0=\mat{s}$, the corresponding canonical,  ``dual'' channel square root iteration is:
\begin{eqnarray}\label{cannonical}
\mat{y}_k &\leftarrow& h_\alpha \left[ \mat{y}_{k-1} \cdot \mat{z}_{k-1} \right] \cdot \mat{y}_{k-1}  \nonumber \\
\mat{z}_k &\leftarrow& \mat{z}_{k-1} \cdot h_\alpha \left[ \mat{y}_{k-1} \cdot \mat{z}_{k-1} \right] \; ,
\end{eqnarray}
converging as ${\mat{y}}_k \rightarrow \mat{s}^{1/2}$, ${\mat{z}}_k \rightarrow \mat{s}^{-1/2}$ and
${\mat{x}}_k \rightarrow {\mat{I}}$, with eigenvalues aggregated towards 1 by
the NS map $h_\alpha[\mat{x}]=\frac{\sqrt{\alpha}}{2} \left(3-\alpha \mat{x} \right)$ \cite{Higham08,higham2005}.
As in the case of sign iteration, this canonical iteration was shown by Higham, Mackey, Mackey and
Tisseur \cite{higham2005} to remain strongly bounded in the super-linear regime, by idempotent {Fr\'{e}chet} derivatives about the fixed point
$\left(\mat{s}^{1/2},\mat{s}^{-1/2}\right)$, in the direction $\left(
\delta \mat{y}_{k-1} , \delta \mat{z}_{k-1} \right)$:
\begin{eqnarray}
\delta \mat{y}_k &=& \frac{1}{2} \delta \mat{y}_{k-1} - \frac{1}{2} \mat{s}^{1/2} \cdot \delta \mat{z}_{k-1} \cdot \mat{s}^{1/2} \\
\delta \mat{z}_k &=& \frac{1}{2} \delta \mat{z}_{k-1} - \frac{1}{2} \mat{s}^{-1/2} \cdot \delta \mat{y}_{k-1} \cdot \mat{s}^{-1/2} \;.
\end{eqnarray}
In addition to the dual channel instance, we also consider the ``single'' channel version of square root iteration,
\begin{eqnarray}\label{single}
\mat{z}_k &\leftarrow& \mat{z}_{k-1} \cdot h_\alpha \left[ \mat{x}_{k-1} \right] \; , \nonumber \\
\mat{x}_k &\leftarrow&  \mat{z}^\xpose_{k} \cdot \mat{s} \cdot \mat{z}_{k} \; .
\end{eqnarray}

\section{Implementation}\label{implementation}

\subsection{Programming}

In our experimental research, issue driven implementations of $\tt SpAMM$ have been developed,
including a Haskell version (formal functional programming) \cite{spammh},
a fine grained ($\tt 4\times 4$) {{\tt single-precision}} assembly coded version (scalar performance) \cite{Bock2013} and
a task parallel version in {\tt C++}, {\tt OpenMP 3.0} and {\tt Charm++} (strong scaling) \cite{BockCK14}.
In the current contribution, informal functional programming in {\tt Fortran08} was used,
with the goal of generic simplicity and  mathematical agility.

In our implementation, allocation functions  instantiate or reuse sub-matrices in {downward} recursion,
and accumulate  decorations (flops, bounding boxes, non-zeros, norms, initialization flags {{\em etc.}})  in
backwards recursion, up the stack. Optional, $\tt ifdef$'d features include the first order {Fr\'{e}chet}
analyses, outlined in Section \ref{stability} and using $\tt MATMUL$,
as well as sparse $\tt VTK$ output for visualization of the $ijk$ product volumes, shown in Section \ref{locality}.

Precision is determined by the block dimension $N_b$, the primary threshold $\tau$
controlling error in the $\mat{z}$ and the $\mat{x}$ channels, and by the tighter (sensitive) threshold
$\tau_s$ for the $\mat{y}$ channel.  Unless stated otherwise, we take $N_b=16$ and  $\tau_s  \sim {\tt .01}\times \tau$.
Finally, reported calculations were carried out in double precision using the {\tt GCC/gfortran 4.8.1} compiler.

\subsection{Mapping}\label{map}
The NS logistic map for the square root iteration is $h_\alpha[\mat{x}]=\frac{\sqrt{\alpha}}{2} \left(3-\alpha \mat{x} \right)$,
with the initial rate of convergence controlled by $h'_\alpha$ and the smallest eigenvalue, $x_0$.
Various schemes for controlling the values $\alpha$ towards convergence include methods by
Pan and Schreiber \cite{Pan1991}, and more recently, Jie and Chen \cite{chen2014}, who demonstrated  $\tt 2 \times$
acceleration for very ill-conditioned problems with their {continuous} scaling approach.

In addition to scaling of the NS logistic, we introduce a stabilizing map that accounts for eigenvalues
tossed out of bounds by $\ot$. This stabilization is the transformation $[0,1] \rightarrow [0+\varepsilon, 1-\varepsilon]$
(shift and scale), carried out prior to application of the logistic.

The most important aspect of these scaling and stabilization maps is to turn them off towards convergence.  Conventional methods
often compute a lowest eigenvalue to monitor convergence \cite{Pan1991,chen2014}, but this may be too expensive for
ill-conditioned problems.  Alternatively, we monitor convergence simply with the relative trace error,
$t_k = \left( n - {\rm tr} \, \widetilde{\mat{x}}_k \right) n^{-1}$.
Then, sigmoidal functions damp scaling to unity,
\begin{equation}
\alpha(t) = {\tt 1.} + {\tt 1.85}  \times \left( 1+e^{-{\tt 50.}(t-{\tt .35})}  \right)^{-1}  \, ,
\end{equation}
and the stability parameter to zero,
\begin{equation}
\varepsilon(t) = {\tt .1} \times \left( 1+e^{-{\tt 75.}(t-{\tt .30})}  \right)^{-1} \, .
\end{equation}
These empirical damping functions are used throughout.

\subsection{Data} \label{data}

Data for numerical experiments include problems from electronic structure and structural engineering.
Electronic structure matrices {were obtained} from the non-orthogonal metric (overlap matrix) of the generalized eigenproblem,
encountered in local support with Gaussian-Type Atomic-Orbitals (GTAOs) \cite{helgaker2008molecular}.
A sequence of nanotube matrices, $36 \times \rightarrow 128 \times$ the (3,3) unit cell,
were generated with {\tt TubeGen} \cite{tubegen} and a modified STO-2G \cite{Schuchardt2007} basis,
with an added diffuse (flat) Gaussian $sp$-shell and exponents $\zeta_{10^{10}}={\tt .06918}$ and $\zeta_{10^{11}}={\tt .05934}$,
corresponding to the {condition numbers} $\kappa(\mat{s})=10^{10}$ and $\kappa(\mat{s})=10^{11}$
respectively\footnote{In this case, $\kappa$ is double exponential with decreasing $\zeta$.}.
We also constructed a sequence of matrices from periodic water cubes, in increments of $100$,
with coordinates obtained using the {\tt gromacs} utility {\tt gmx solvate -box} \cite{gmx30}
and the triple-$\zeta$ 6-311G** GTAO basis \cite{Schuchardt2007}.
While less ill-conditioned than the nano-tube sequence,
$\kappa(\mat{s})\sim 10^5$, the water cube matrices manifest a different metric locality due to dimensionality.
Also, we experiment with the {\tt bcsstk14} structural engineering matrix for the Roof of the Omni Coliseum \cite{BCSSTK14}.

\section{Error Flow} \label{errorflow} 

\subsection{Stability}\label{stability}


Stability in the square root iteration is determined  by the differential
\begin{equation} \label{firstorderdual}
\delta \mat{x}_k = \,  { \mat{x}}_{\delta \widehat{\mat{y}}_{k-1}}  \, {\scriptstyle \times} \, \delta y_{k-1}
                 \, + \,  { \mat{x}}_{\delta \widehat{\mat{z}}_{k-1}}  \, {\scriptstyle \times} \, \delta z_{k-1}  \, + {\cal O}(\tau^2) \, ,
\end{equation}
which must remain bounded below one to avoid divergence.   The corresponding {Fr\'{e}chet} derivatives are
\begin{equation}
  \mat{x}_{\delta \widehat{ \mat{y}}_{k-1}}
= \lim_{\tau \rightarrow 0} \frac{ \mat{x} \left( \mat{y}_{k-1} + \tau \delta \widehat{\mat{y}}_{k-1} ,  {\mat{z}}_{k-1}  \right)
                                     -\mat{x}_k    }{\tau}
 \end{equation}
and
 \begin{equation}
 \mat{x}_{\delta \widehat{ \mat{z}}_{k-1}} = \lim_{\tau \rightarrow 0}
\frac{ \mat{x} \left( {\mat{y}}_{k-1} , \mat{z}_{k-1} +\tau  \delta \widehat{\mat{z}}_{k-1} \right) - \mat{x}_k   }{\tau}  \, ,
 \end{equation}
along unit directions of the previous errors $\delta \widehat{\mat{y}}_{k-1}$ and $\delta \widehat{\mat{z}}_{k-1}$, by an amount
determined by the displacements $\delta y_{k-1} = \lVert \delta \mat{y}_{k-1} \rVert$  and  $\delta z_{k-1}=\lVert \delta \mat{z}_{k-1} \rVert$.
In the single instance, we have simply:
\begin{equation} \label{firstordersingle}
\delta \mat{x}_k = \,  { \mat{x}}_{\delta \widehat{\mat{z}}_{k-1}}  \, {\scriptstyle \times} \, \delta z_{k-1}  \, + {\cal O}(\tau^2)  \, .
\end{equation}

This formulation makes plain changes about the resolvent, separating orientational effects for
derivatives of the unit direction, set mostly by the underlying exact linear algebra, from
changes to error displacements, which involve both the action of derivatives on previous errors,  as well as
current $\tt SpAMM$ occlusion errors.  In the following sections, we develop this form of the error.  Then,
in Section \ref{mayss}, we show interesting behaviors of these derivatives at the edge of stability.

\subsection{{Fr\'{e}chet} derivatives}

In the dual instance, {Fr\'{e}chet} derivatives occurring in Eq.~(\ref{firstorderdual}) are:


\begin{multline}\label{dxdy}
 \mat{x}_{\delta \widehat{ \mat{z}}_{k-1}} =  {\mat{y}}_{k-1} \cdot  h'_\alpha \delta \widehat{ \mat{z}}_{k-1} \cdot  \mat{y}_{k-1}  \cdot \mat{z}_{k}
 + \mat{y}_k \cdot  \delta \widehat{\mat{z}}_{k-1} \cdot   h_\alpha \left[ \mat{x}_{k-1} \right] \\
 +  \mat{y}_{k} \cdot  \mat{z}_{k-1} \cdot {\mat{y}}_{k-1} \cdot h'_\alpha \delta \widehat{\mat{z}}_{k-1} \, ,
\end{multline}

and
\begin{multline}\label{dxdz}
  \mat{x}_{\delta \widehat{ \mat{y}}_{k-1}} = h_\alpha \left[ \mat{x}_{k-1} \right]  \cdot \delta \widehat{\mat{y}}_{k-1} \cdot \mat{z}_{k}
+  h'_\alpha  \delta \widehat{\mat{y}}_{k-1} \cdot \mat{z}_{k-1} \cdot  \mat{y}_{k-1} \cdot  \mat{z}_{k} \\
 + \mat{y}_{k} \cdot \mat{z}_{k-1} \cdot h'_\alpha \delta \widehat{\mat{y}}_{k-1} \cdot \mat{z}_{k-1}  \, .
\end{multline}

Closer to the fixed point orbit,  $\mat{y}_k \cdot \mat{z}_{k-1} \rightarrow \mat{I}$, $\mat{y}_{k-1} \cdot \mat{z}_{k} \rightarrow \mat{I}$,
$h_\alpha \left[ \mat{x}_{k} \right] \rightarrow \mat{I}$ and $h'_\alpha \rightarrow - \frac{1}{2}$ \cite{higham2005}.  Then,

\begin{equation} \label{yorbit}
 \mat{x}_{\delta \widehat{ \mat{y}}_{k-1}} \rightarrow \delta \widehat{\mat{y}}_{k-1} \cdot \left( \mat{z}_k-\mat{z}_{k-1} \right)
\end{equation}
and
\begin{equation} \label{zorbit}
 \mat{x}_{\delta \widehat{ \mat{z}}_{k-1}} \rightarrow \left( \mat{y}_k-\mat{y}_{k-1} \right) \cdot \delta \widehat{\mat{z}}_{k-1} .
\end{equation}
Likewise, in the single channel instance:
\begin{equation}
 \mat{x}_{\widehat{\mat{z}}_{k-1}} \rightarrow  \left(  \mat{z}_{k} - \mat{z}_{k-1} \right)^\xpose \cdot \mat{s} \cdot \delta \widehat{\mat{z}}_{k-1}
+  \, \delta \widehat{\mat{z}}^\xpose_{k-1} \cdot  \mat{s}  \cdot \left(  \mat{z}_{k} - \mat{z}_{k-1} \right)  \, .
\end{equation}
About the fixed point then, error flow in the $\mat{y}$ and the $\mat{z}$ channels is tightly quenched,
corresponding to $\mat{x}_{\delta \widehat{\mat{x}}_{k-1}} \rightarrow  \mat{I}$  and identity iteration \cite{higham2005}.

\begin{figure}[h]
\includegraphics[width=5.in]{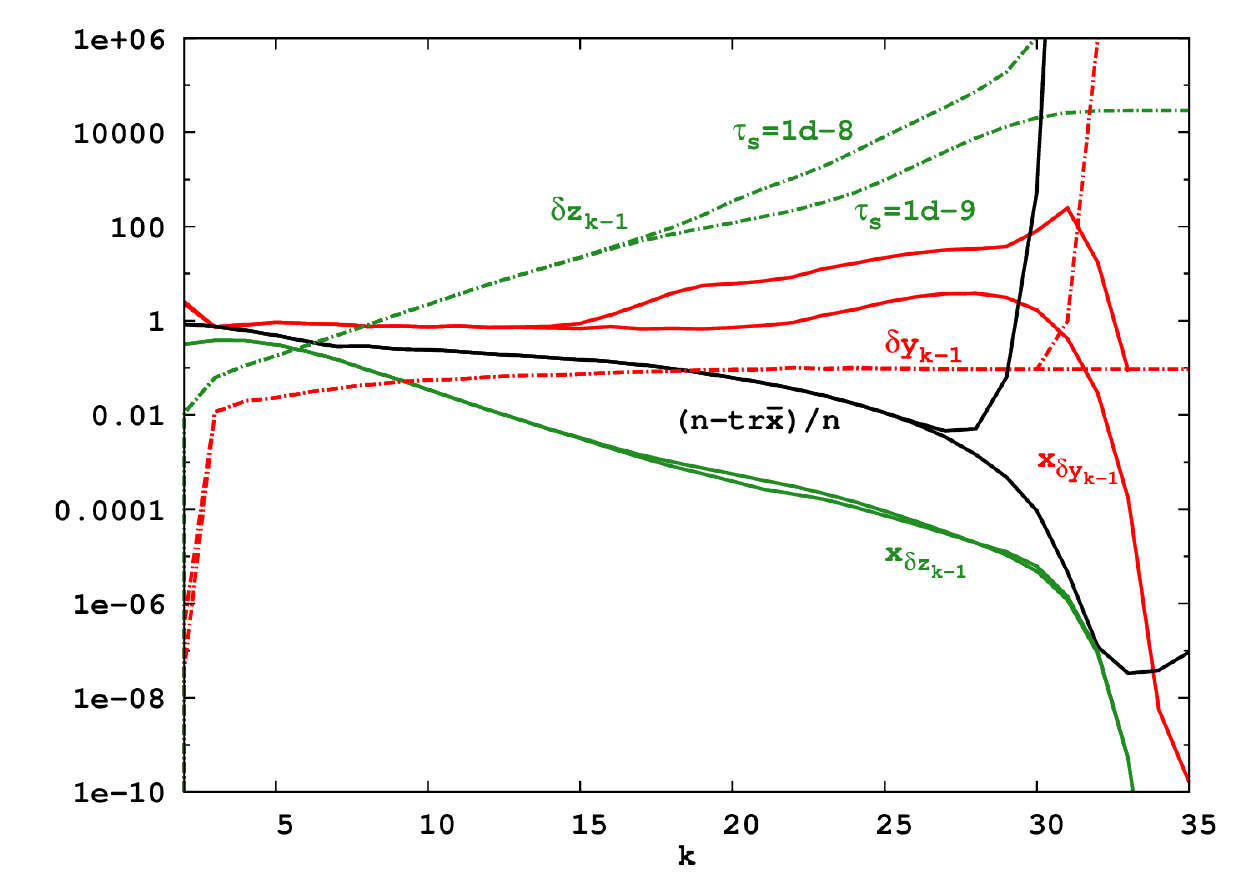}
\caption{Trace error and $\lVert \cdot \rVert$ of derivatives and displacements for the unscaled dual iteration.
Derivatives are full lines, whilst displacements for $\tau_s=\{10^{-8}, 10^{-9}\}$
are dashed lines.  The trace error is a full black line. } \label{flow_noscale_dual}
\end{figure}

\begin{figure}[h]
\includegraphics[width=5.in]{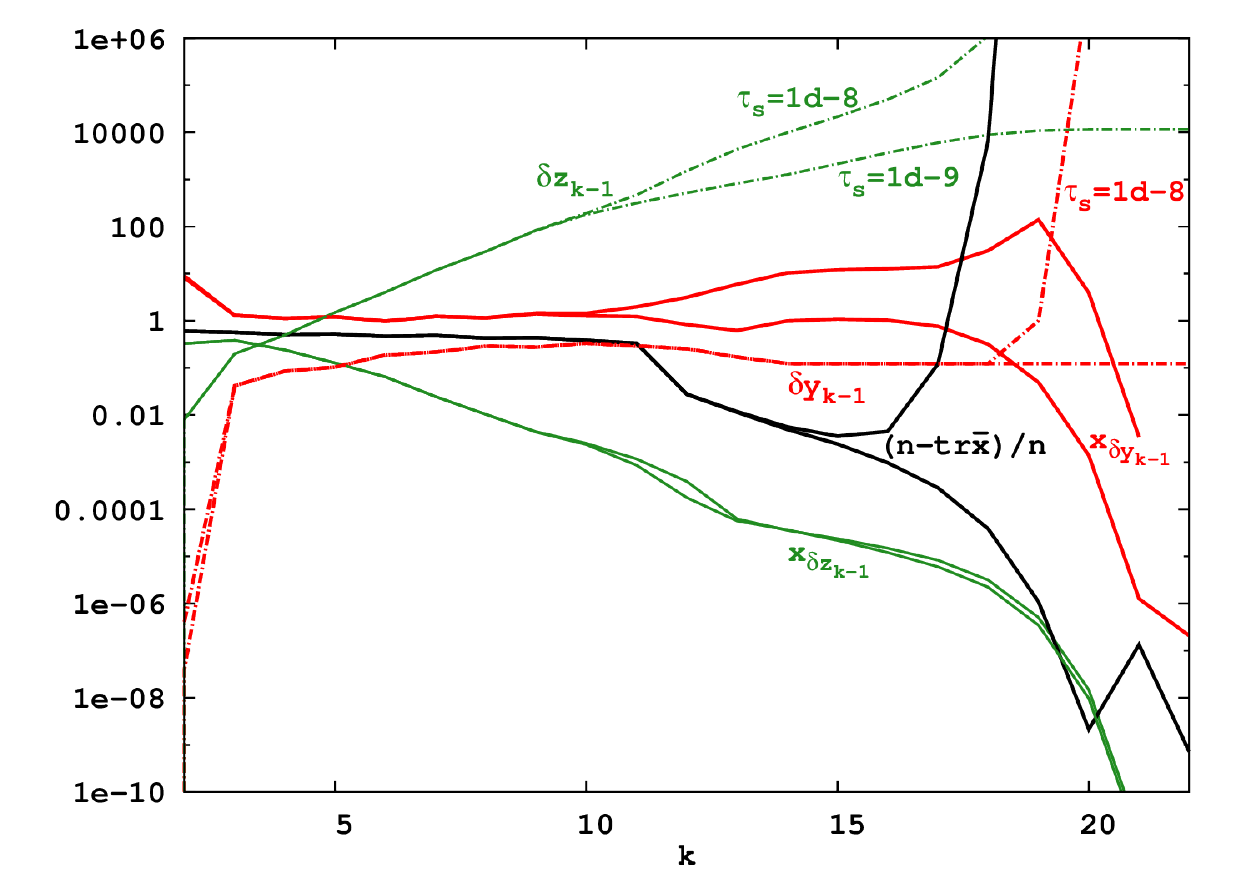}
\caption{Trace error and $\lVert \cdot \rVert$ of derivatives and displacements for the scaled dual iteration.
Derivatives are full lines, whilst displacements for $\tau_s=\{10^{-8}, 10^{-9}\}$
are dashed lines.  The trace error is a full black line. }\label{flow_scaled_dual}
\end{figure}

\begin{figure}[h]
\includegraphics[width=5.in]{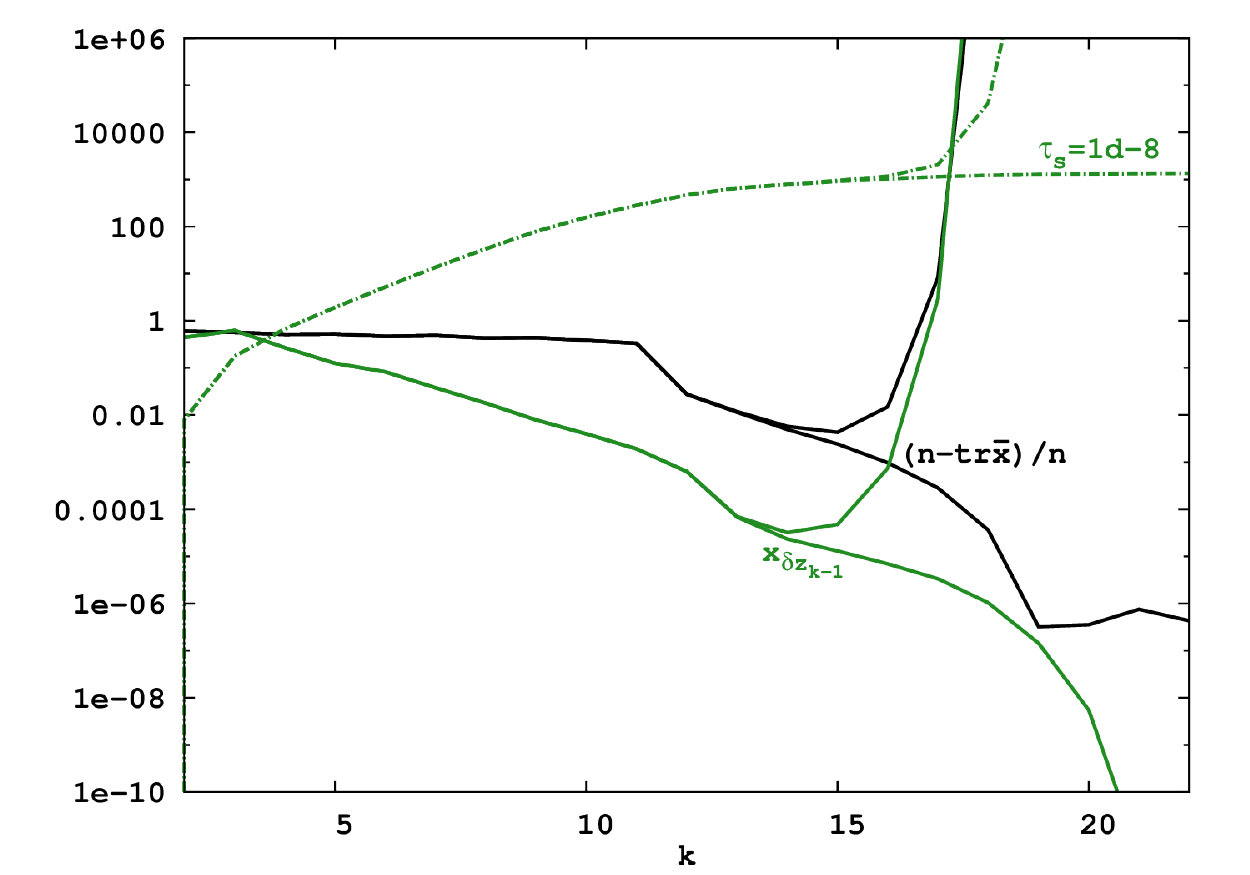}
\caption{Trace error and $\lVert \cdot \rVert$ of derivatives and displacements for the scaled single iteration.
Derivatives are full lines, whilst displacements for $\tau_s=\{10^{-7}, 10^{-8}\}$
are dashed lines.  The trace error is a full black line. }\label{flow_scaled_stab}
\end{figure}

\subsection{Displacements}
At each step, the accumulation of previous errors in addition to the $\tt SpAMM$ occlusion error move the
approximate iteration away from the unperturbed reference, here the {\tt {double-precision}} iteration of arrays with {\tt MATMUL}.

Including the $\tt SpAMM$ error in the $\widetilde{\mat{z}}_{k-1}$ update we have:
\begin{equation} \label{widetildez}
 \widetilde{\mat{z}}_{k-1} =  \widetilde{\mat{z}}_{k-2}  \, \ot \, h_\alpha[\widetilde{\mat{x}}_{k-2}]
= \Delta^{\widetilde{\mat{z}}_{k-2} \cdot h_\alpha \left[ \widetilde{\mat{x}}_{k-2}\right]}_\tau
+ \widetilde{\mat{z}}_{k-2} \cdot h_\alpha\left[ \widetilde{\mat{x}}_{k-2}\right] \, .
\end{equation}
Then, with $ h_\alpha \left[ \widetilde{\mat{x}}_{k-2} \right]
=  h_\alpha \left[ \mat{x}_{k-2} \right] +  h'_\alpha  \delta \mat{x}_{k-2}$, and taking $\mat{z}_{k-1}$ from both sides,
\begin{equation}
 \delta {\mat{z}}_{k-1} =\Delta^{\widetilde{\mat{z}}_{k-2} \cdot h_\alpha \left[ \widetilde{\mat{x}}_{k-2}\right]}_\tau
 +\delta \mat{z}_{k-2} \cdot h_\alpha \left[\widetilde{\mat{x}}_{k-2} \right]
+ \mat{z}_{k-2} \cdot h'_\alpha \delta \mat{x}_{k-2}  \, ,
\end{equation}
which is bounded by
 \begin{multline}\label{zdispalcementbound}
  \delta {z}_{k-1} <
 \lVert \mat{z}_{k-2} \rVert \left( \tau  \,  n^2  \, \lVert h_\alpha \left[\widetilde{\mat{x}}_{k-2} \right]  \rVert
 + h'_\alpha  \delta y_{k-2} \lVert z_{k-2} \rVert \right)  \\
 + \delta {z}_{k-2} \left( \lVert h_\alpha \left[\widetilde{\mat{x}}_{k-2}  \right] \rVert  + \lVert y_{k-2} \rVert \right) .
 \end{multline}

In Eq.~(\ref{zdispalcementbound}),  the term $h'_\alpha  \delta y_{k-2} { \lVert \mat{z}_{k-2} \rVert }^2$ is volatile, tending towards
$\delta y_{k-2} \, \kappa(\mat{s})/2$.  Because of this sensitivity, and because the $\mat{y}$ product channel maintains fidelity
of the starting eigen-basis, we single out this ``sensitive'' product for a higher level of precision; $\tau_s \ll \tau$.

In the single instance, the $\mat{y}$ channel is implicit in the first product involving $\mat{s}$, which can be from the left or the
right.  In this work, the most accurate product in the single instance is rightmost.

\subsection{Most approximate {yet} still stable}\label{mayss}

The potential to compute fast and effective preconditioners with $\tt SpAMM$ is
determined by the most approximate yet still stable (MAYSS) iteration, a challenge for increasingly ill-conditioned problems.
Illustrative experiments were carried out on the $\kappa(\mat{s})=10^{10}$ nanotube examples described in Section~\ref{data}.
We picked  $\tau ={\tt .001}$ and $N_b = 32$.  Then, we looked at stability with respect to the tighter $\tau_s$ threshold:
In Fig.~\ref{flow_noscale_dual}, unscaled results for the dual instance are shown.  In Fig.~\ref{flow_scaled_dual}, scaled results for the
dual instance are given, and in Fig.~\ref{flow_scaled_stab} we show results for the scaled single instance.

In the dual instances, Figs.~\ref{flow_noscale_dual} \& \ref{flow_scaled_dual}, the bifurcating orientational components of the
error manage to avoid the numerical catastrophe, with $\mat{x}_{\widehat{\mat{z}}_{k-1}}$ in solid green converging strongly,
and $\mat{x}_{\widehat{\mat{y}}_{k-1}}$ in solid red, with an above unity drift driving divergence of the displacements (dashed lines).
On the other hand, bifurcation in the single
instance ($\tau_s=10^{-7}$) finds the orientational component of the error, $\mat{x}_{\delta \widehat{\mat{z}}_{k-1}}$, diverging
well ahead of the displacement $\delta z_{k-1}$.

In these problems, values of $\tau$ near the MAYSS bifurcation do not lead to a reduced complexity; instead, near total fill of the product
is observed towards convergence, even for the largest 128$\times$ unit cell nano-tube.  Also, scaling as reported in Section \ref{implementation} reclaims
about $\sim 2/3$ of the available $2 \times$ acceleration possible at this level off ill-conditioning \cite{chen2014}, but
dramatically enhances fill-in of the metric tensor, via the multiplicative effect of $h'_\alpha$ in Eq.~(\ref{zdispalcementbound}).
In addition to scaling, the single instance also results in a much larger volumetric fill-in, involving extended, delocalized
error flows in the orbit.

Our interpretation of these results is that despite a similar overall convergence behavior and error control,
the tensor volumes accessed by the two instances is very different, due to the magnitude of norms entering the $\tt SpAMM$ kernel;  in the dual instance
$y^{\rm dual}_k=h[\mat{x_{k-1}}] \ot \mat{y}_{k-1}$ is well behaved , while $y^{\rm single}_k=\mat{z}^\xpose_{k} \ot \mat{s}$ encumbers
large norms associated with the broad spectral resolution, leading to extended delocalization of the metric tensor.


\begin{figure}[h]
\includegraphics[width=5.0in]{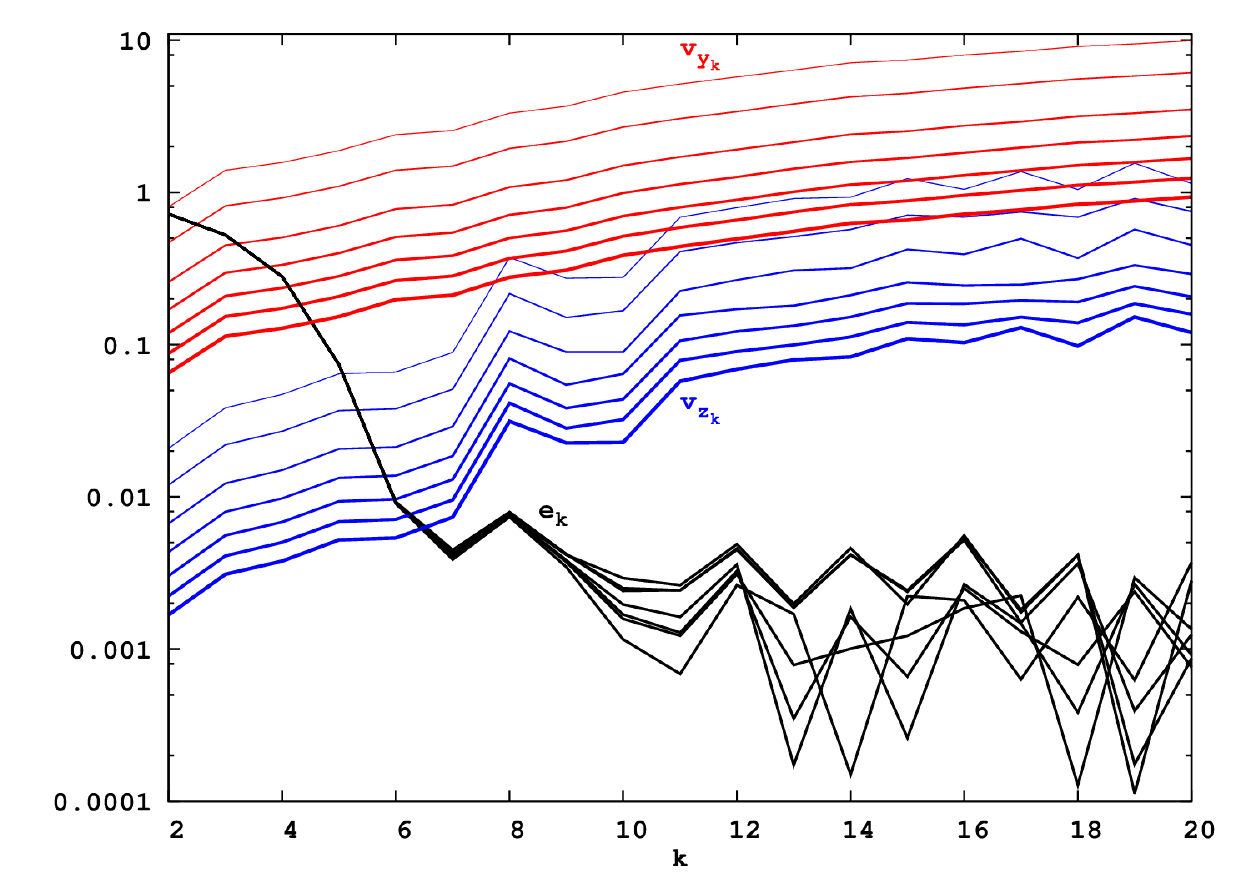}
\caption{
Culled volumes in the thin slice, single instance approximation of $\mat{s}^{-1/2}_{\tau_0 \mu_0}$
for the $\kappa(\mat{s})=10^{10}$ nanotube sequence (line width increasing with system size).
With $\mu_0=\tt .1$ it was only possible to achieve stability
down to $\tau_0=10^{-2} \;  \&  \; \tau_s=10^{-4}$.  Shown are
$\rm v_{\widetilde{\mat{z}}_k} \equiv \left( {\rm vol}_{ \widetilde{\mat{z}}_{k-1}\ot h[\widetilde{\mat{x}}_{k-1}] } \right) \times 100\% / n^3$ (blue) and
${\rm v}_{\widetilde{\mat{y}}_k} \equiv \left( {\rm vol}_{\mat{s}  \ots  \widetilde{\mat{z}}_{k}} \right) \times 100\% / n^3$ (red).
Also shown is the trace error, ${\rm t}_{k} = \left( n-{\rm tr}\,\widetilde{\mat{x}}_k \right)/n$ (black).}\label{regularized_stab}
\end{figure}

\begin{figure}[t]
\includegraphics[width=5.0in]{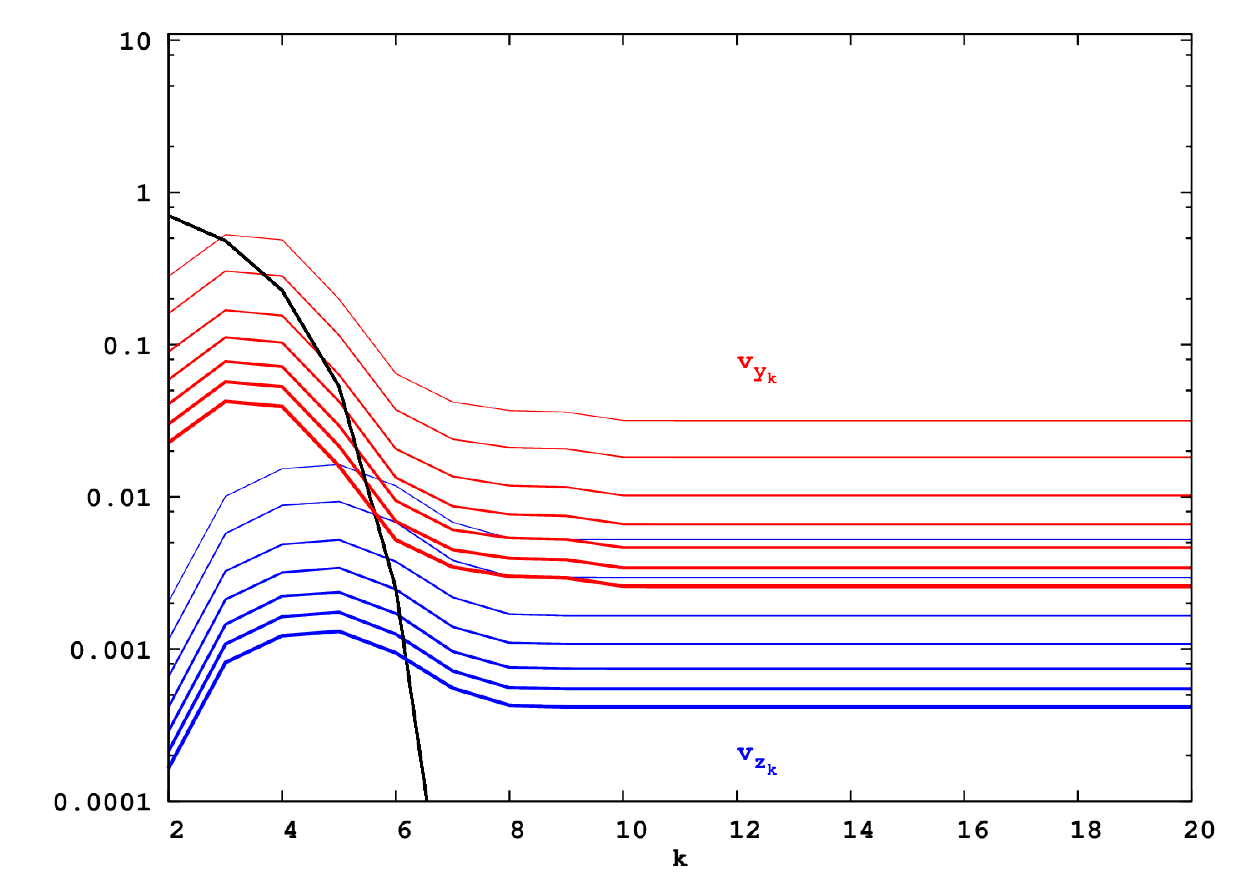}
\caption{
Culled volumes in the thin slice, dual instance approximation of $\mat{s}^{-1/2}_{\tau_0 \mu_0}$
for the $\kappa(\mat{s})=10^{10}$ nanotube sequence (line width increasing with system size).
The slice is $\mu_0={\tt .1}, \tau_0={\tt .1} \;  \&  \; \tau_s={\tt .001}$ thin.
Shown are
${\rm v}_{\widetilde{\mat{y}}_k}= \left( {\rm vol}_{  h[\widetilde{\mat{x}}_{k-1}] \ots \widetilde{\mat{y}}_k }  \right) \times 100\% / n^3$ (blue) and
${\rm v}_{\widetilde{\mat{z}}_k}= \left( {\rm vol}_{\widetilde{\mat{z}}_{k-1} \ot  h[\widetilde{\mat{x}}_{k-1}]} \right) \times 100\% / n^3$ (red).
Also shown is the trace error, ${\rm t}_{k} = \left( n-{\rm tr}\,\widetilde{\mat{x}}_k \right)/n$ (black),
which rapidly approaches $10^{-11}$ (not shown).}\label{regularized_dual}
\end{figure}

\section{Regularization}\label{regularization}

Even for the most approximate {yet still} stable approximations (MAYSS), our nanotube calculations lead to
delocalized products that are not tightly bound by Eq.~(\ref{bound}), even for very large {128$\times$ unit cell} systems.
And while similarly ill-conditioned problems may achieve substantial compression with just the MAYSS approximation,
as shown later in Fig.~\ref{Lensing4}, the $\tt SpAMM$ approximation cannot generally yield a fast method in cases of
severe ill-conditioning.

A systematic way to reduce these effects is through Tikhonov regularization \cite{neumaier1998solving,sarra2014}.
Regularization involves a small level shift of the eigenvalues,  $\mat{s}_\mu \leftarrow \mat{s}+\mu \mat{I}$, altering the
condition number of the shifted matrix to  $\kappa( \mat{s}_\mu) = \frac{\sqrt{s^2_{n-1} + \mu^2}}{\sqrt{s^2_0+\mu^2}}$ \cite{sarra2014}.

Achieving substantial acceleration with severe ill-conditioning  may require a large level shift however,
producing inverse factors of little practical use.  One approach to recover a more accurate inverse
factor is Riley's method based on Taylor's series \cite{sarra2014};
\begin{equation}
\mat{s}^{-1/2} = \mat{s}^{-1/2}_{\mu} \cdot \left( \mat{I}+\frac{\mu}{2} \mat{s}^{-1}_{\mu}
                                                   +\frac{3 \mu^2}{8} \mat{s}^{-2}_{\mu} + \dots
   \right) \; .
\end{equation}
For severely ill-conditioned problems and large level shifts, this expansion may converge very slowly.
Also, adding powers of the full inverse may not be computationally effective.

\subsection{Product representation}

We introduce an alternative representation of the regularized inverse factor;
\begin{equation} \label{spammsandwich}
\mat{s}^{-1/2} \equiv \bigotimes_{\substack{\tau=\tau_0 \\ \mu=\mu_0   } } {\left|\, \tau\, \mu \, ; \, \scriptstyle{\mat{s}^{-1/2}}  \right>}  \, ,
\end{equation}
which is a telescoping product of preconditioned ``slices''
starting with a most-approximate-yet-still-effective-by-one-order (MAYEBOO) preconditioner,
$\mat{s}^{-1/2}_{\tau_0 \mu_0} \equiv {\left|\, \tau_0\, \mu_0 \, ; \, \scriptstyle{\mat{s}^{-1/2}}\right>}$ 
Braket notation marks the potential for asymmetries in the intermediate representation.
This sandwich of generic, thinly sliced {\tt SpAMM} products allows to construct a nested scoping on precision via $\tau$, and in the
effective condition number controlled by $\mu$.


\subsection{Effective by one order}

We look again at the $\kappa(\mat{s})=10^{10}$ nanotube series described in Section~\ref{data},
this time with extreme regularization, $\mu_0={\tt .1}$, and at a finer granularity, $N_b=8$.
Culled $\mat{y}_k$ and $\mat{z}_k$ volumes (as percentage of the total work) for $36 - 128 \, \times$ the (3,3) unit cell
are shown for the MAYEBOO approximation in Fig.~\ref{regularized_stab} for the single instance,
and in Fig.~\ref{regularized_dual} for the dual instance.

The behavior of these implementations  is very different; in the single
instance, a stable  iteration could not be found at precision $\tau_0=\tt .1$.  Stability could
only be found at ${\tt .01}$, and that with a poorly contained trace error and cull-volumes
that continue to inflate past convergence, with a conspicuous $\sqrt{k}$-like dependence.
This behavior results from the broad resolution of spectral powers
$\widetilde{\mat{y}}_k \rightarrow  \mat{s}^{-1/2}_{\tau_0 \mu_0} \oto \mat{s}_{\mu_0}$,
with corresponding large metric fields that are poorly bound by Eq.~(\ref{bound}).

On the other hand, dual iteration volumes collapse rapidly with fast trapping of the trace error,
as $\widetilde{\mat{y}}_k\rightarrow  \mat{I}_{\tau_0 \mu_0} \oto \mat{s}^{1/2}_{\tau_0 \mu_0}$
and $\widetilde{\mat{z}}_k \rightarrow  \mat{s}^{-1/2}_{\tau_0 \mu_0} \oto \mat{I}_{\tau_0 \mu_0}$,
and with Eq.~(\ref{bound}) tightening to
\begin{equation}\label{boundY}
\Delta^{\mat{I}_{\tau_0 \mu_0} \cdot \mat{s}^{1/2}_{\tau_0 \mu_0}} <  \tau n \lVert \mat{s}^{1/2}_{\tau_0 \mu_0} \rVert
\end{equation}
and
\begin{equation}\label{boundZ}
\Delta^{ \mat{s}^{-1/2}_{\tau_0 \mu_0}\cdot \mat{I}_{\tau_0 \mu_0}}  <  \tau n \lVert \mat{s}^{-1/2}_{\tau_0 \mu_0} \rVert \, .
\end{equation}
This contraction to the plane diagonal is compressive, leading to computational complexities that should approach quadtree copy in place.


\begin{figure}[h]
\fbox{\includegraphics[width=3.85cm,keepaspectratio=true,
                       trim={.5cm 2.3cm 2.cm 1.cm},clip]
                       {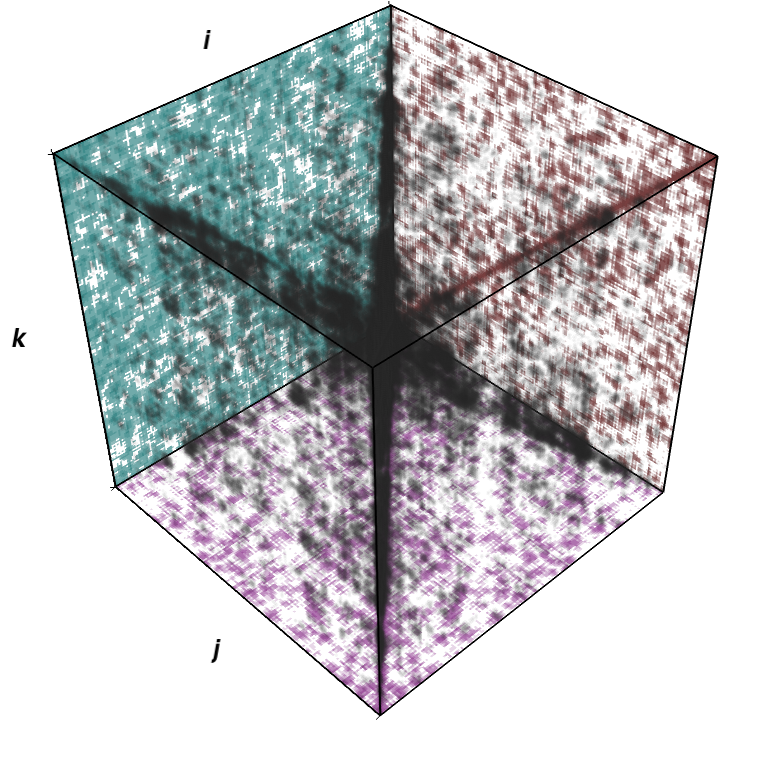}}
\fbox{\includegraphics[width=3.85cm,keepaspectratio=true,
                        trim={0.5cm 2.3cm 2.cm 1.cm},clip]
                        {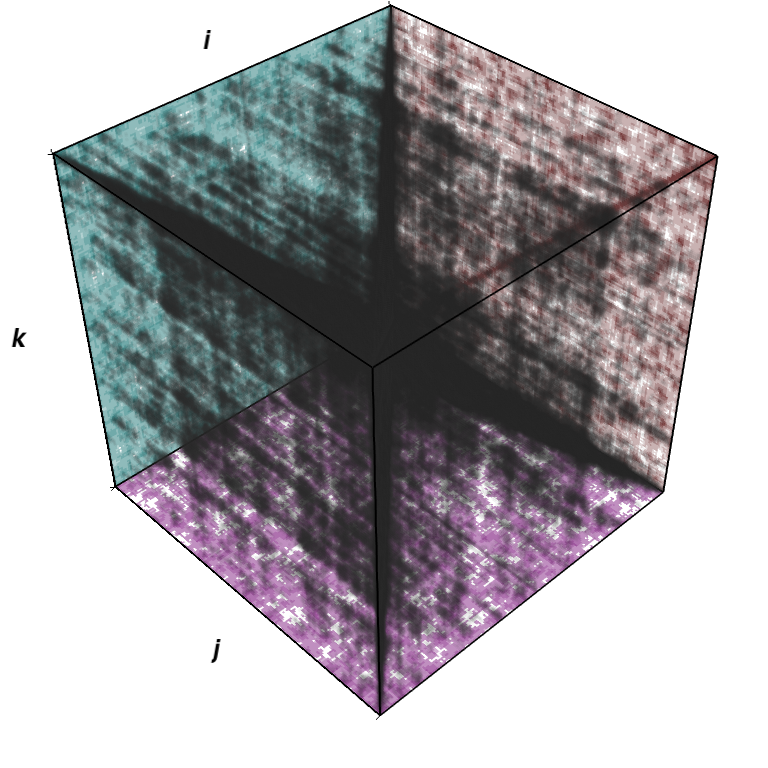}}
\fbox{\includegraphics[width=3.85cm,keepaspectratio=true,
                        trim={0.5cm 2.3cm 2.cm 1.cm},clip]
                        {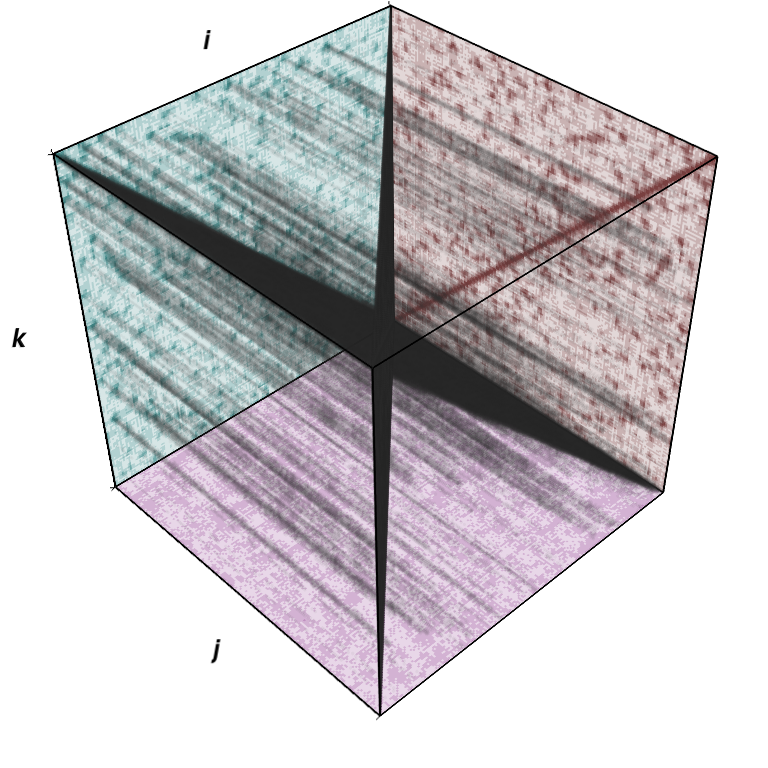}}

\fbox{\includegraphics[width=3.85cm,keepaspectratio=true,
                        trim={0.5cm 2.3cm 2.cm 1.cm},clip]
                        {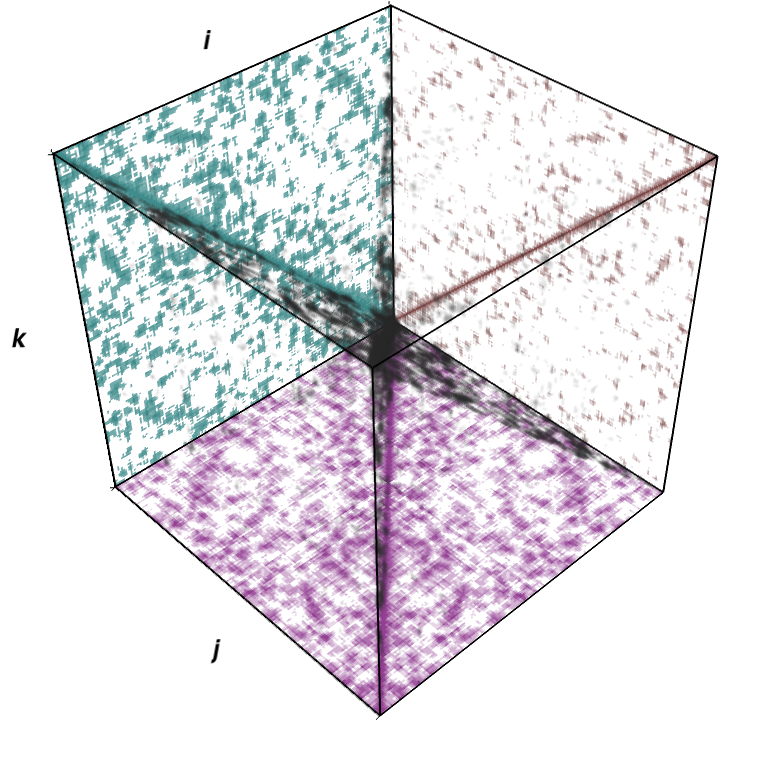}}
\fbox{\includegraphics[width=3.85cm,keepaspectratio=true,
                        trim={0.5cm 2.3cm 2.cm 1.cm},clip]
                        {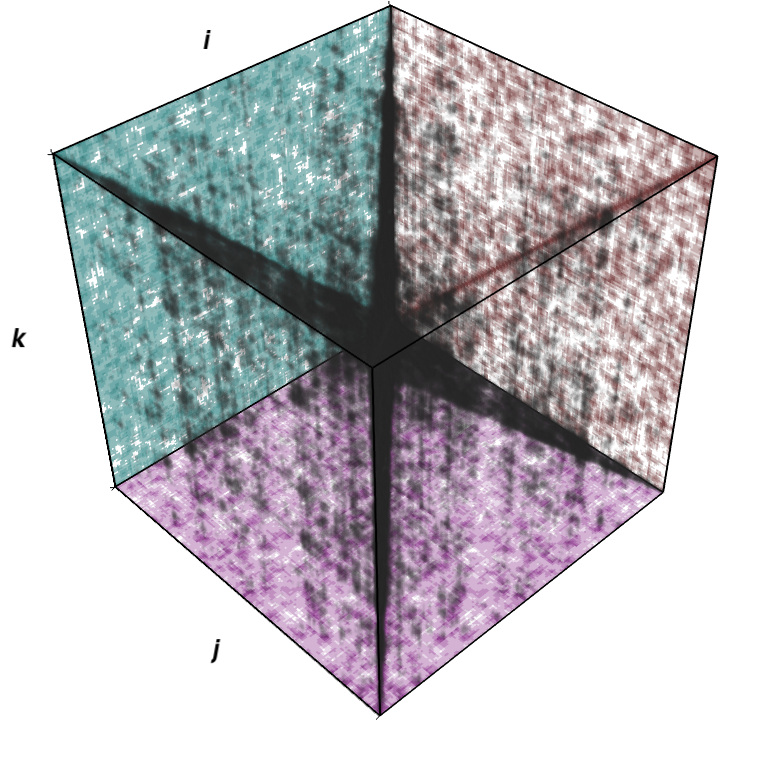}}
\fbox{\includegraphics[width=3.85cm,keepaspectratio=true,
                        trim={0.5cm 2.3cm 2.cm 1.cm},clip]
                        {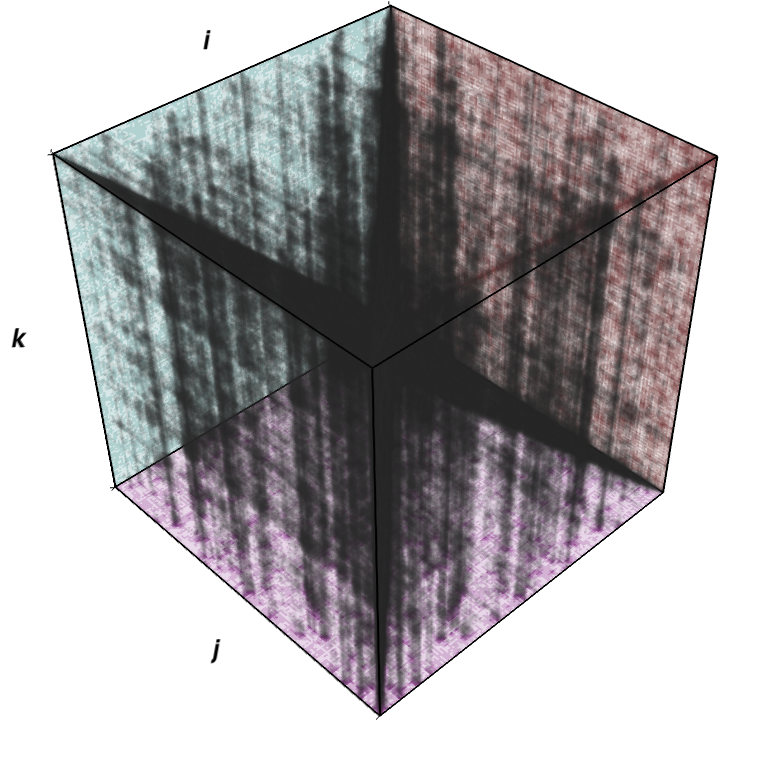}}
\caption{
Product volumes in construction of the unregularized preconditioner
$\left|\tau_0={\tt .001} ; \,\scriptstyle{\mat{s}^{-1/2}} \right>$, with
dual instance square root iteration, and for the 6-311G** metric of 100 periodic water molecules.
At top its  $\mat{y}_k=h_\alpha[ \mat{x}_{k-1} ] \ots \mat{y}_{k-1}$
for $k=0,5,\& 17$, while on the bottom we have $\mat{x}_k=  \mat{y}_{k}  \ot \mat{z}_{k}$ for $k=0,5, \& 17$.
Maroon is $\mat{a}$, purple is $\mat{b}$, green is $\mat{c}$,  and black is the volume ${\rm vol}_{a \ot b}$
in the product $\mat{c}=\mat{a} \ot \mat{b}$.}\label{Lensing3}
\end{figure}

\subsection{Iterative regularization}

We now sketch an iterative approach to constructing the product representation, Eq.~(\ref{spammsandwich}).
In the dual instance, it is possible to obtain a first MAYEBOO approximation $\mat{s}^{-1/2}_{\tau_0={\tt .1}, \mu_0={\tt .1}}$,
which improves the condition number by one order of magnitude, with a numerical resolution of approximately one digit.
Then, a next level slice can be found, $\mat{s}^{-1/2}_{\tau_0 \mu_1}$, based on the residual
$\left(\mat{s}^{-1/2}_{\tau_0\mu_0} \right)^\xpose \, \otone \, \left(\mat{s}+\mu_1 \mat{I} \right)
\, \otone \,\mat{s}^{-1/2}_{\tau_0 \mu_0} $, with {\em e.g.}~%
$\mu_1= \tt .01$ and $\tau_1={\tt .01}$.   The product $\mat{s}^{-1/2}_{\tau_0 \mu_1}  \otone \mat{s}^{-1/2}_{\tau_0 \mu_0}$
then improves the condition number by two orders of magnitude, still with a numerical resolution of one digit.
Reflected in the preceding notation, it appears necessary to compute the residual at a higher level of precision,
{\em e.g.} using $\otone$ instead of $\oto$ and with $\tau_1 > \tau_0$.

In this way,  it may be possible to obtain product representation of the inverse square root at a $\tt SpAMM$ resolution that is
potentially far more permissive than otherwise possible,
\begin{equation} \label{product_rep}
\mat{s}^{-1/2}_{\tau_0} = \mat{s}^{-1/2}_{\tau_0 \mu_n} \, \otone \, \mat{s}^{-1/2}_{\tau_0 \mu_{n-1}} \, \otone \, \dots  \,  \mat{s}^{-1/2}_{\tau_0 \mu_0} \, ,
\end{equation}
assuming ${\tt .1} \ge \mu_0 > \mu_1 > \dots$  Likewise, it may also  be possible to obtain the full inverse factor with
increasing numerical resolution as
\begin{equation} \label{product_rep_tau}
\mat{s}^{-1/2} = \mat{s}^{-1/2}_{\tau_m} \, \otpm \,  \mat{s}^{-1/2}_{\tau_{m-1}} \, \otm \dots \, \mat{s}^{-1/2}_{\tau_0} \, ,
\end{equation}
and ${\tt .1} \ge  \tau_0 > \tau_1 > \dots $

Also, with each step a well conditioned generic slice, it may be possible to find a more effective logistic map optimized
for a vanilla distribution of eigenvalues. Finally, relative to the  regularization and precision scoping sketched here,
alternative products are possible that may be far more efficient.  We hope to pursue these efforts in future work.

\section{Locality} \label{locality}

\subsection{Spatial and  metric  locality}

Astrophysical $n$-body algorithms employ range queries over spatial databases to hierarchically discover
and compute approximations that commit only small errors.  Often, these spatial databases are ordered with a
space filling curve (SFC)
\cite{Wise:1984:RMQ:1089389.1089398,springerlink:10.1007/3-540-51084-2_9,
      Samet:1990:DAS:77589,Wise1990,Wise:Ahnentafel,
      Lorton:2006:ABL:1166133.1166134,Samet:2006:DBDS,Adams:2006:SOS,Bock2013,bader13},
which maps points that are close in space to an index where they are also close.
Spatial locality of this type empowers the $\tt SpAMM$ approximation through Block-By-Magnitude
orderings of the sub-space metric.

This metric locality is compressive, but diminished by dimensionality.
In Figure \ref{Lensing3}, we show $\ot$ volumes for square root iteration, corresponding to the Gramian matrix of a small, periodic water box
with the large 6-311G** basis (Section \ref{data}).
In this 3-d periodic case, diminishing Cartesian separations lead to long-skinny {pillae}
and related delocalizations not observed in lower dimensional problems at this modest $\kappa(\mat{s})\sim 10^5$ level of ill-conditioning.
These delocalizations correspond to weakness in Eq.~(\ref{bound}),
and to tighter values of $\tau_s$, required in the MAYSS approximation.
As $n$ becomes large, Cartesian separation will eventually thin these delocalizations
leading to complexity reduction due only to metric decay.

\subsection{Algebraic locality}
\begin{figure}[t]
\fbox{ \includegraphics[width=3.8cm,keepaspectratio=true,
                        trim={3.5cm 3.cm 5.cm 5.cm},clip]
                        {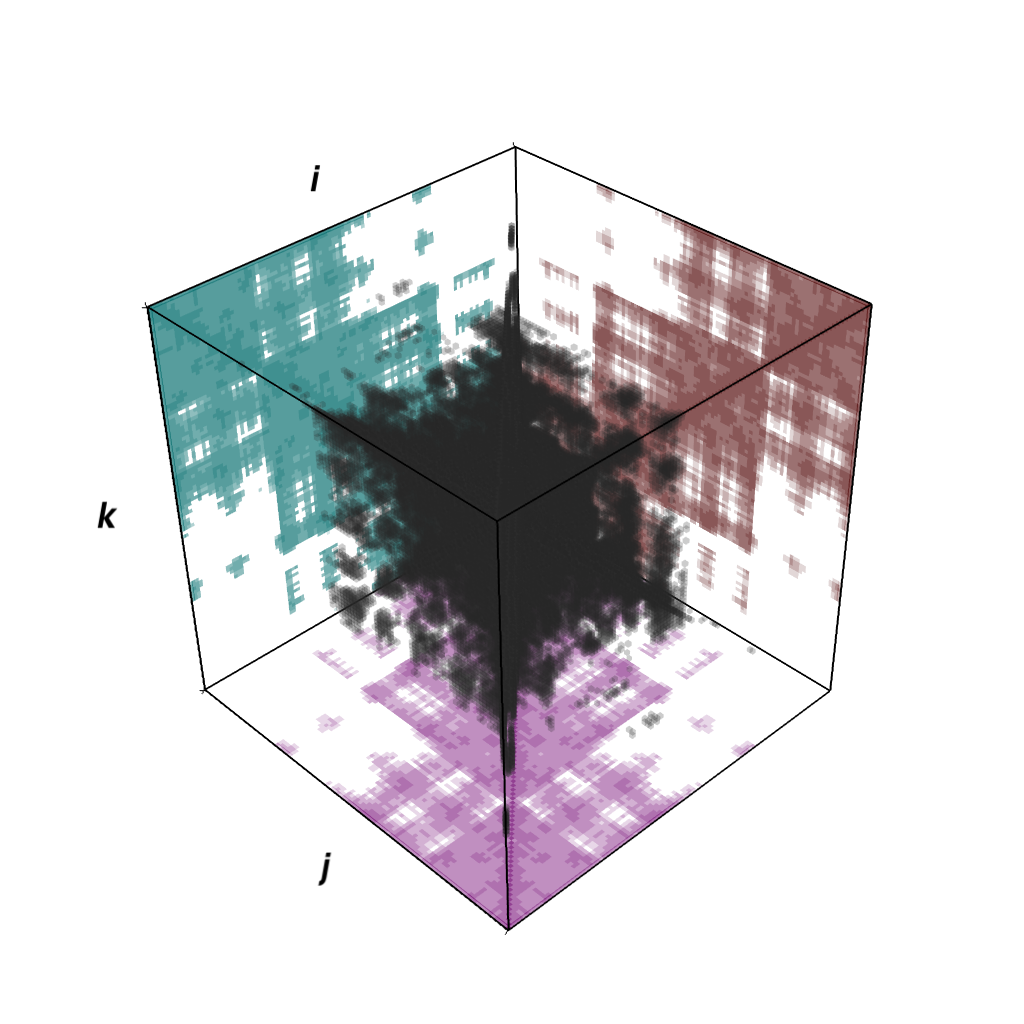}}
\fbox{ \includegraphics[width=3.8cm,keepaspectratio=true,
                        trim={3.5cm 3.cm 5.cm 5.cm},clip]
                        {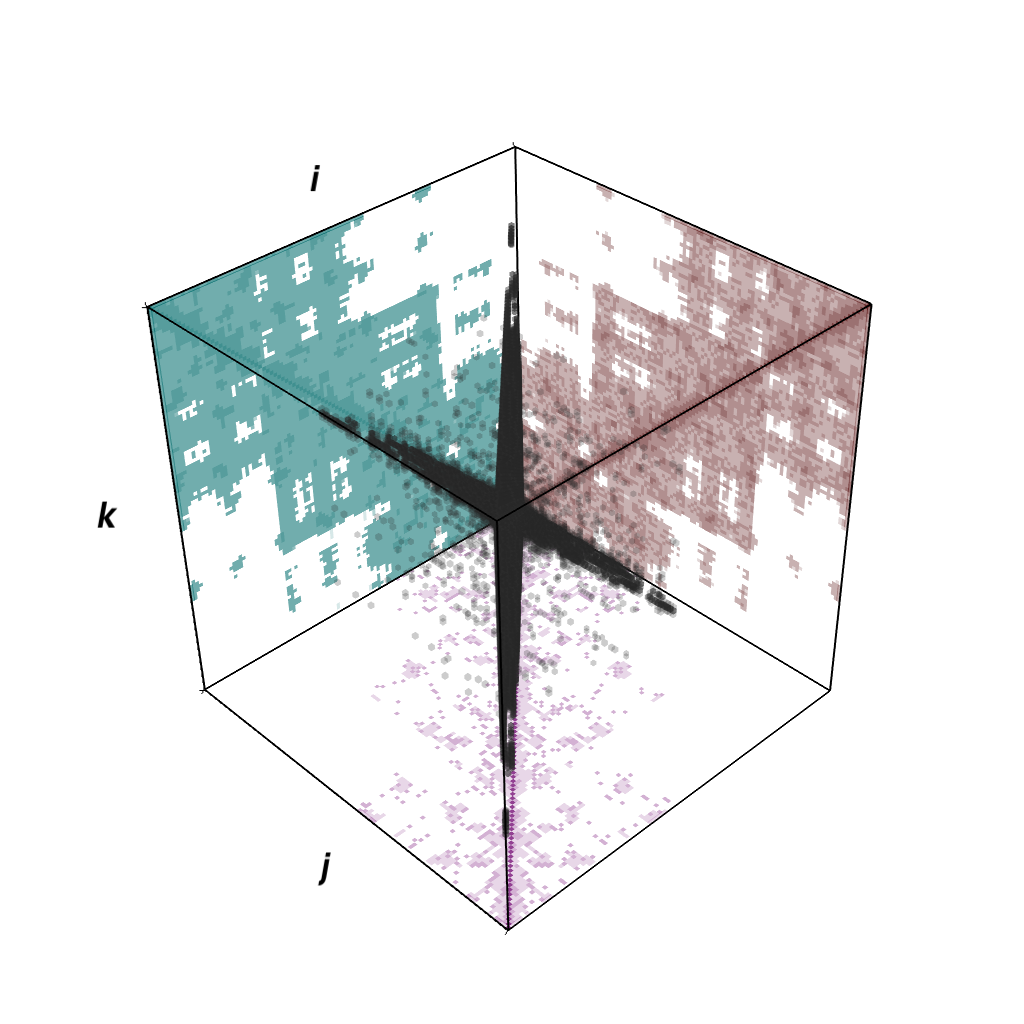}}
\fbox{ \includegraphics[width=3.8cm,keepaspectratio=true,
                        trim={3.5cm 3.cm 5.cm 5.cm},clip]
                        {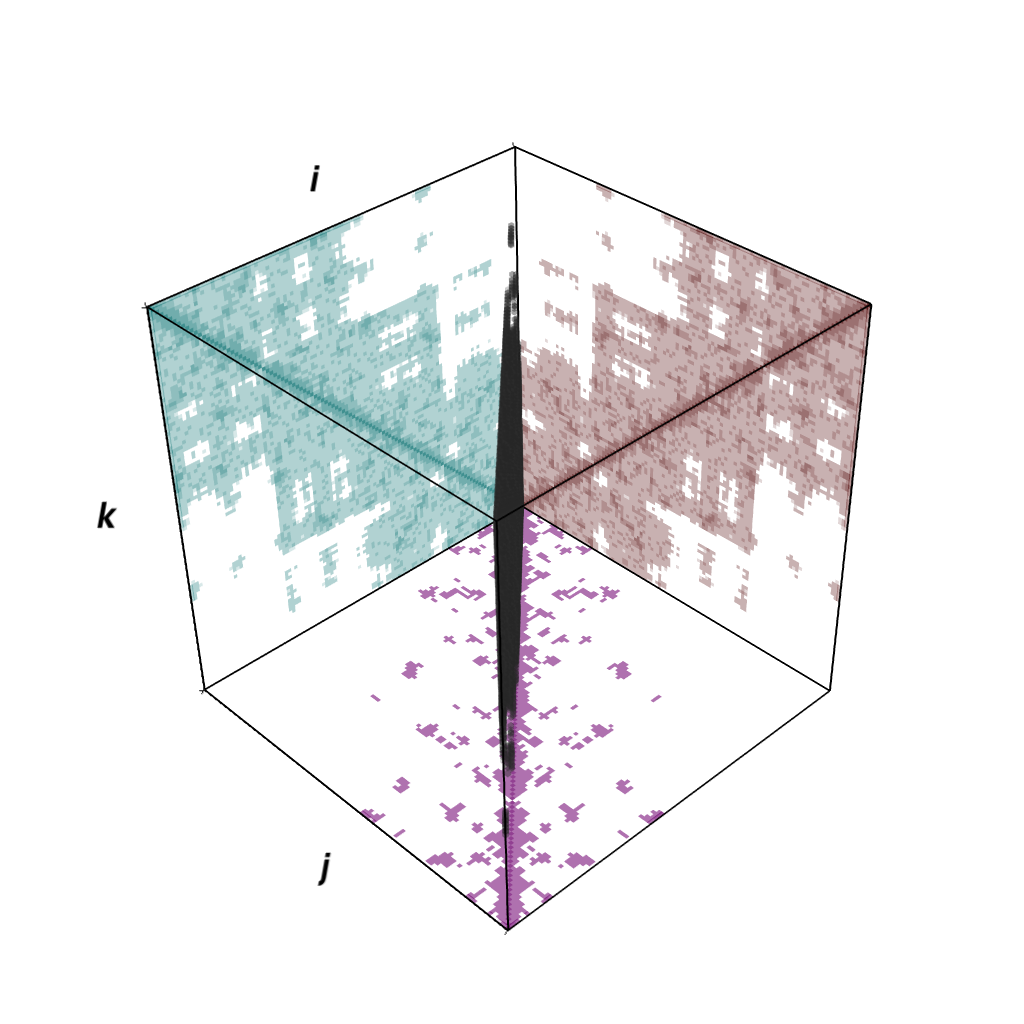}}

\fbox{ \includegraphics[width=3.8cm,keepaspectratio=true,
                        trim={3.5cm 3.cm 5.cm 5.cm},clip]
                        {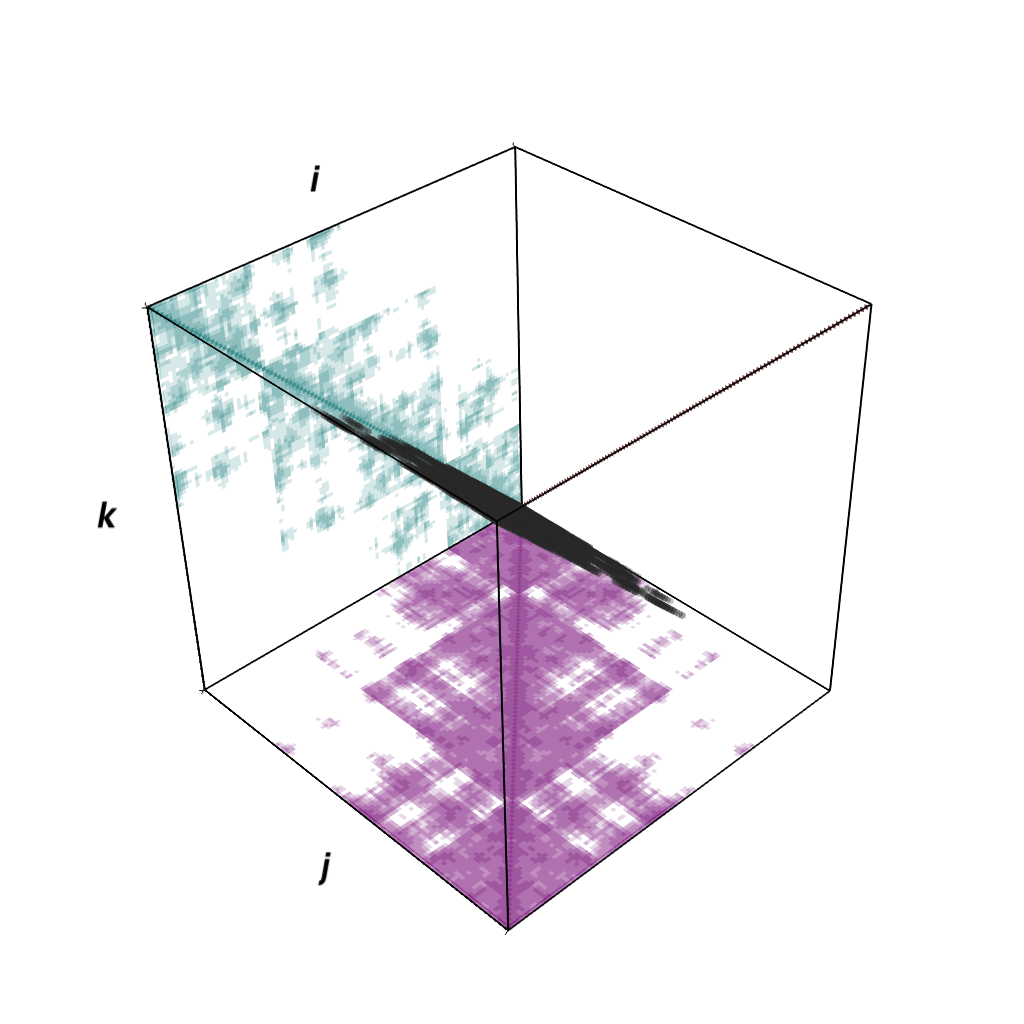}}
\fbox{ \includegraphics[width=3.8cm,keepaspectratio=true,
                        trim={3.5cm 3.cm 5.cm 5.cm},clip]
                        {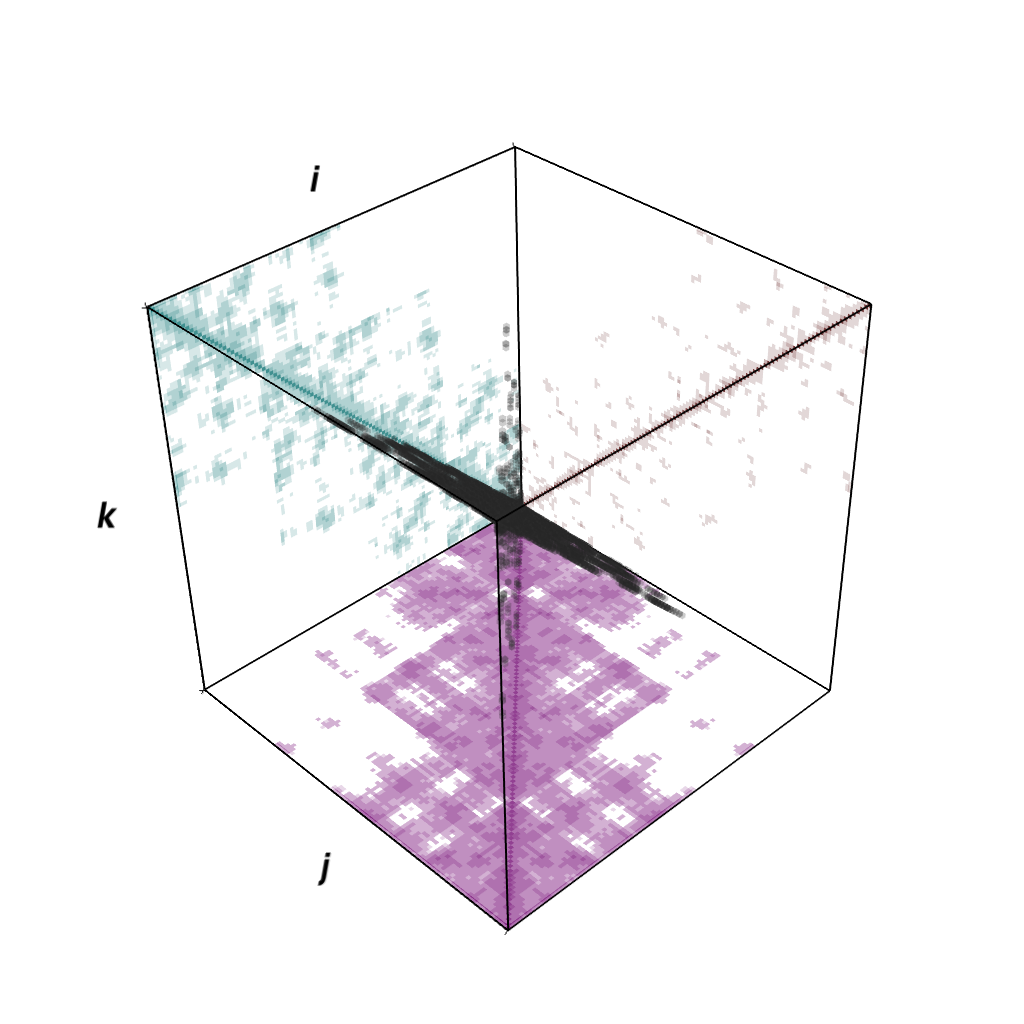}}
\fbox{ \includegraphics[width=3.8cm,keepaspectratio=true,
                        trim={3.5cm 3.cm 5.cm 5.cm},clip]
                        {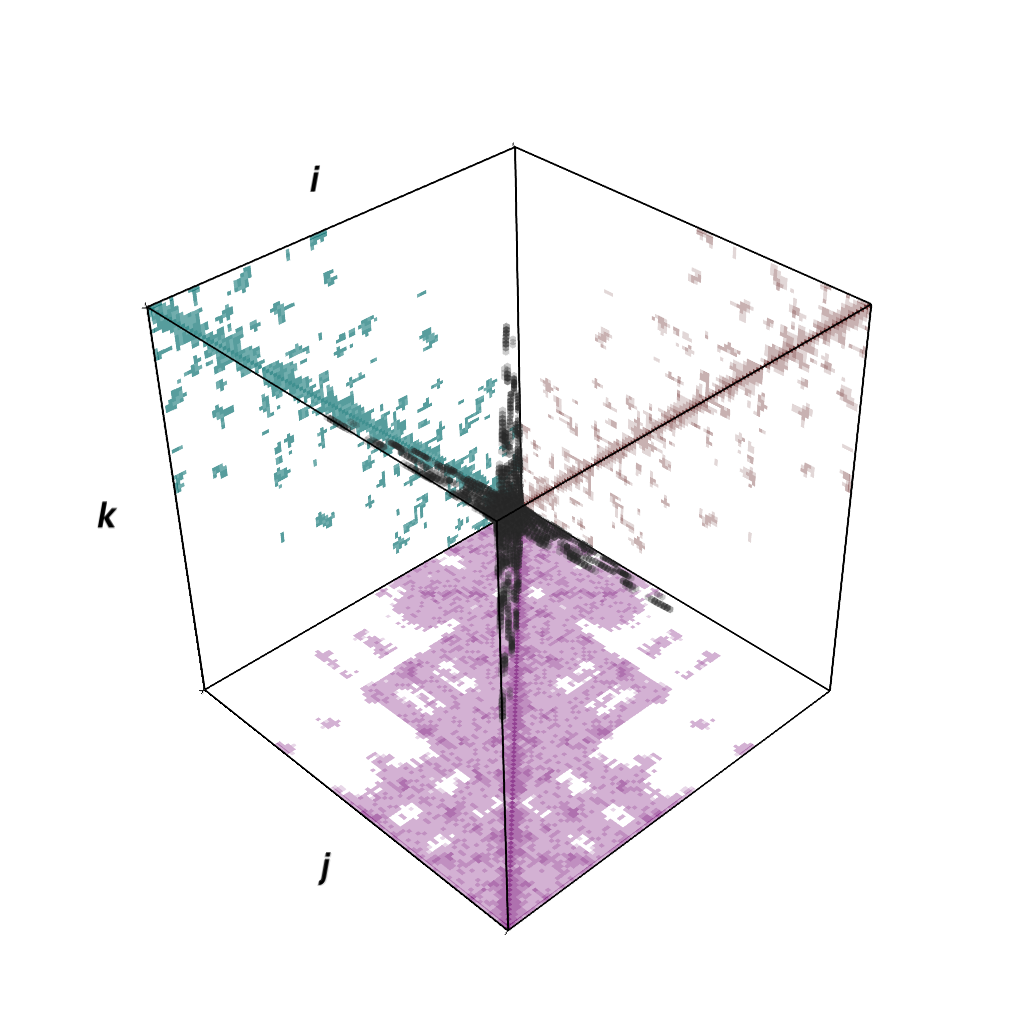}}
\caption{Product volumes in construction of the MAYEBOO preconditioner
$\left|\tau_0=.1,\mu_0=.1\, ; \,\scriptstyle{\mat{s}^{-1/2}} \right>$, with
dual instance square root iteration,  for 8$\times$ $\kappa(s)=10^{11}$ nano-tube.
$\mat{y}_k$ appears wider than $\mat{z}_k$ because it is computed at a higher precision, $\tau_s=\tt .001$,
and because the first multiply involves $\mat{s}^2$.  At top its  $\mat{y}_k=h_\alpha[ \mat{x}_{k-1} ] \ots \mat{y}_{k-1}$
for $k=0,4,\& 16$, while on the bottom we have $\mat{x}_k=  \mat{y}_{k}  \ot \mat{z}_{k}$ for $k=0,2, \& 16$.
Maroon is $\mat{a}$, purple is $\mat{b}$, green is $\mat{c}$,  and black is the volume ${\rm vol}_{a \ot b}$
in the product $\mat{c}=\mat{a} \ot \mat{b}$.}\label{Lensing1}
\end{figure}

\begin{figure}[t]
\fbox{ \includegraphics[width=3.8cm,keepaspectratio=true,
                       trim={0.cm 2.3cm 2.cm 1.cm},clip]
                       {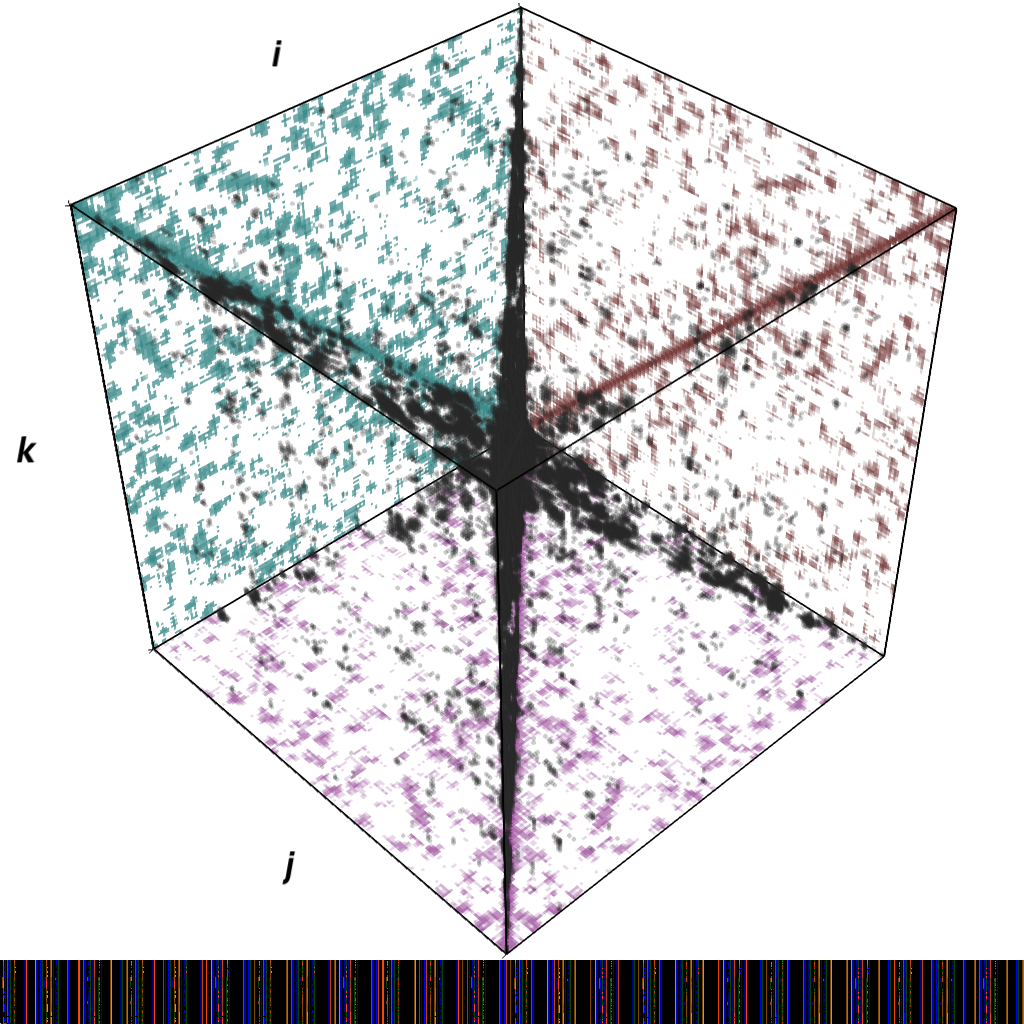}}
\fbox{ \includegraphics[width=3.8cm,keepaspectratio=true,
                        trim={0.cm 2.3cm 2.cm 1.cm},clip]
                        {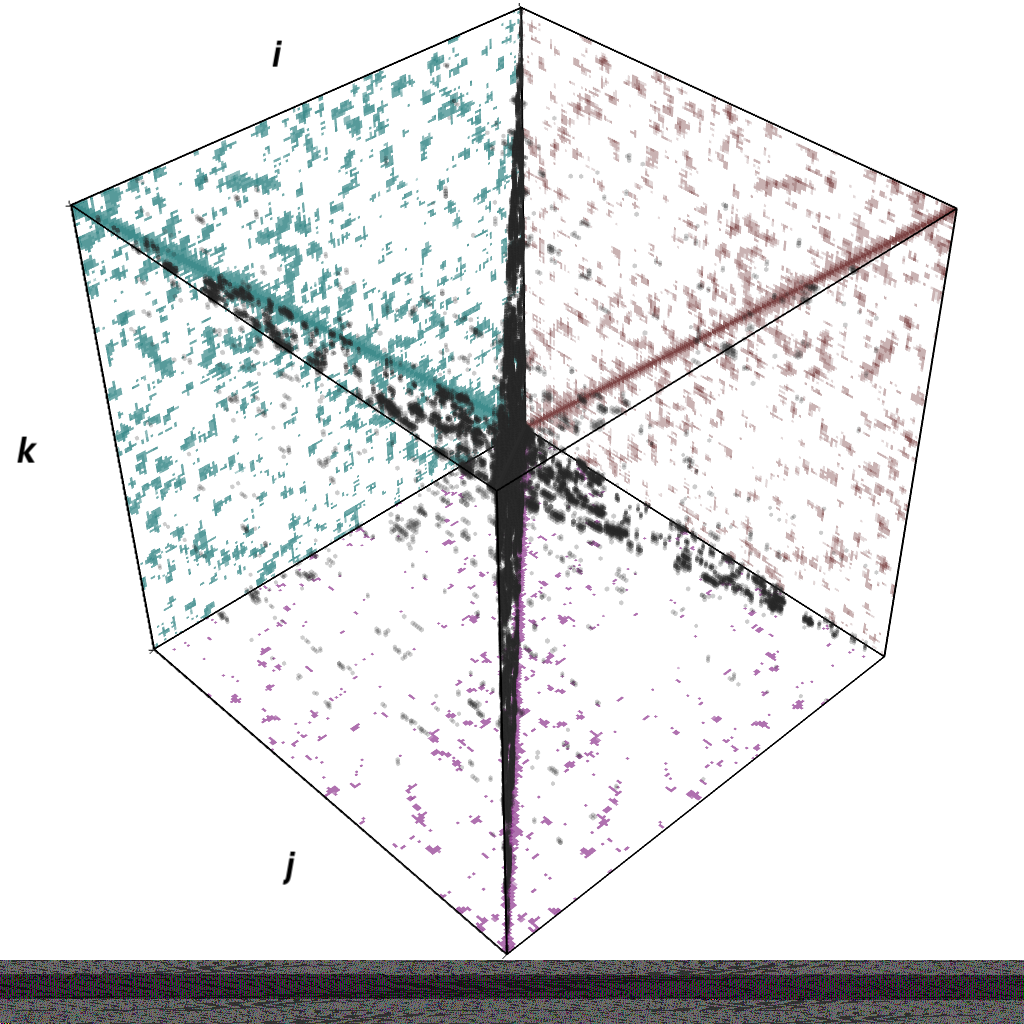}}
\fbox{ \includegraphics[width=3.8cm,keepaspectratio=true,
                        trim={0.cm 2.3cm 2.cm 1.cm},clip]
                        {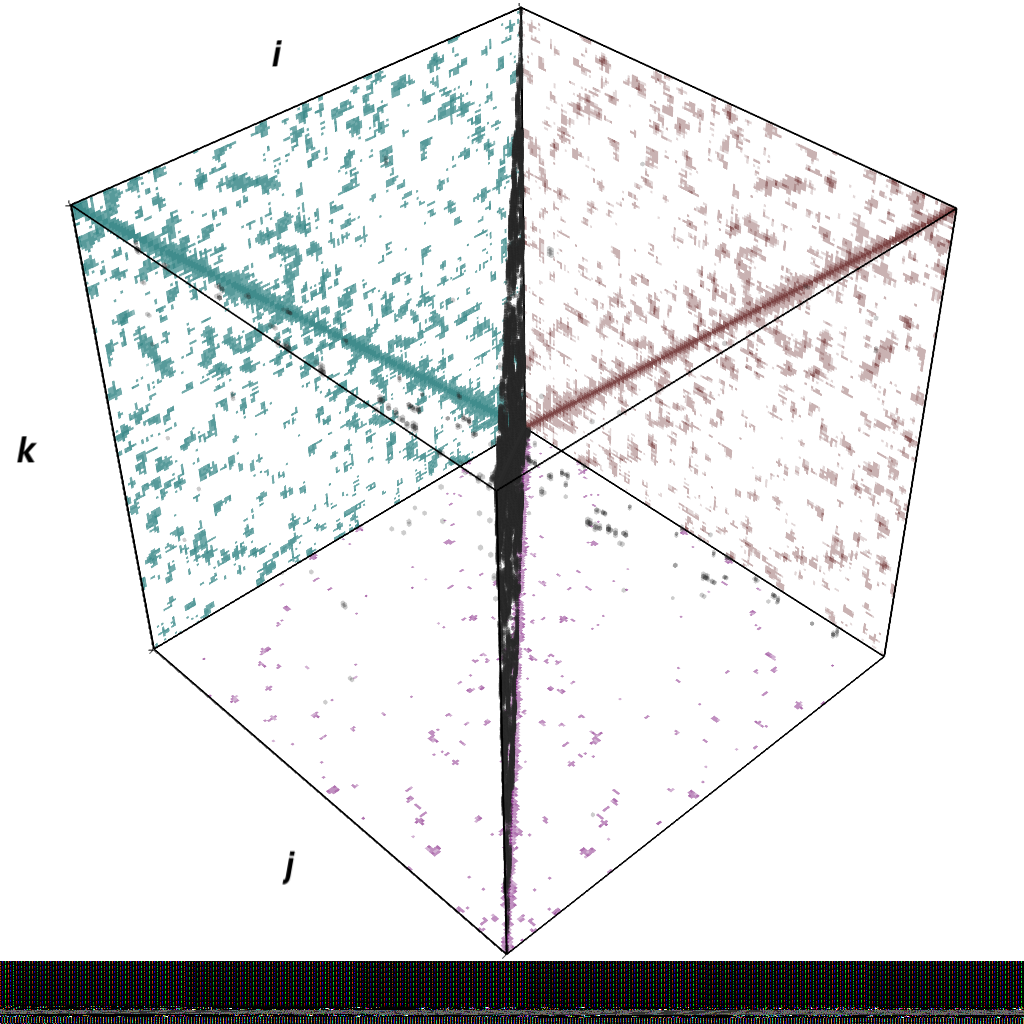}}

\fbox{ \includegraphics[width=3.8cm,keepaspectratio=true,
                        trim={0.cm 2.3cm 2.cm 1.cm},clip]
                        {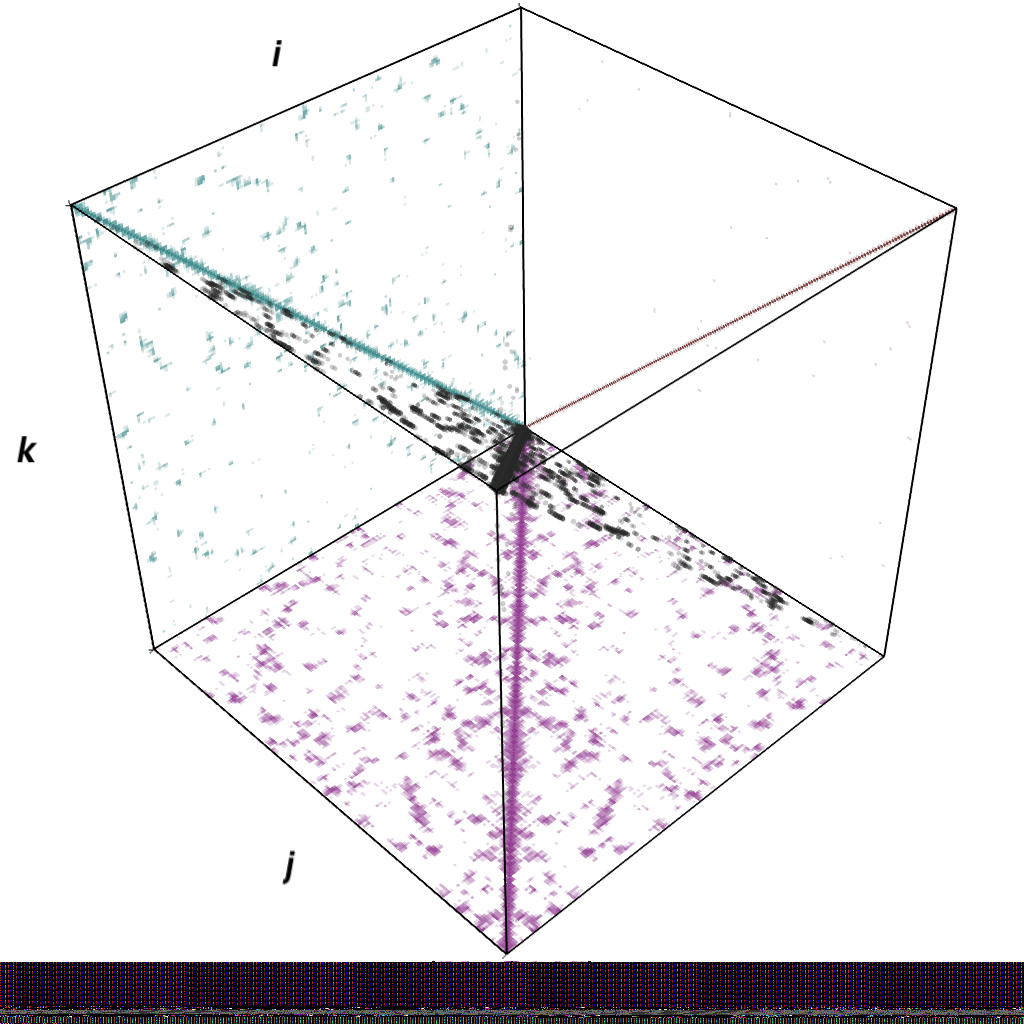}}
\fbox{ \includegraphics[width=3.8cm,keepaspectratio=true,
                        trim={0.cm 2.3cm 2.cm 1.cm},clip]
                        {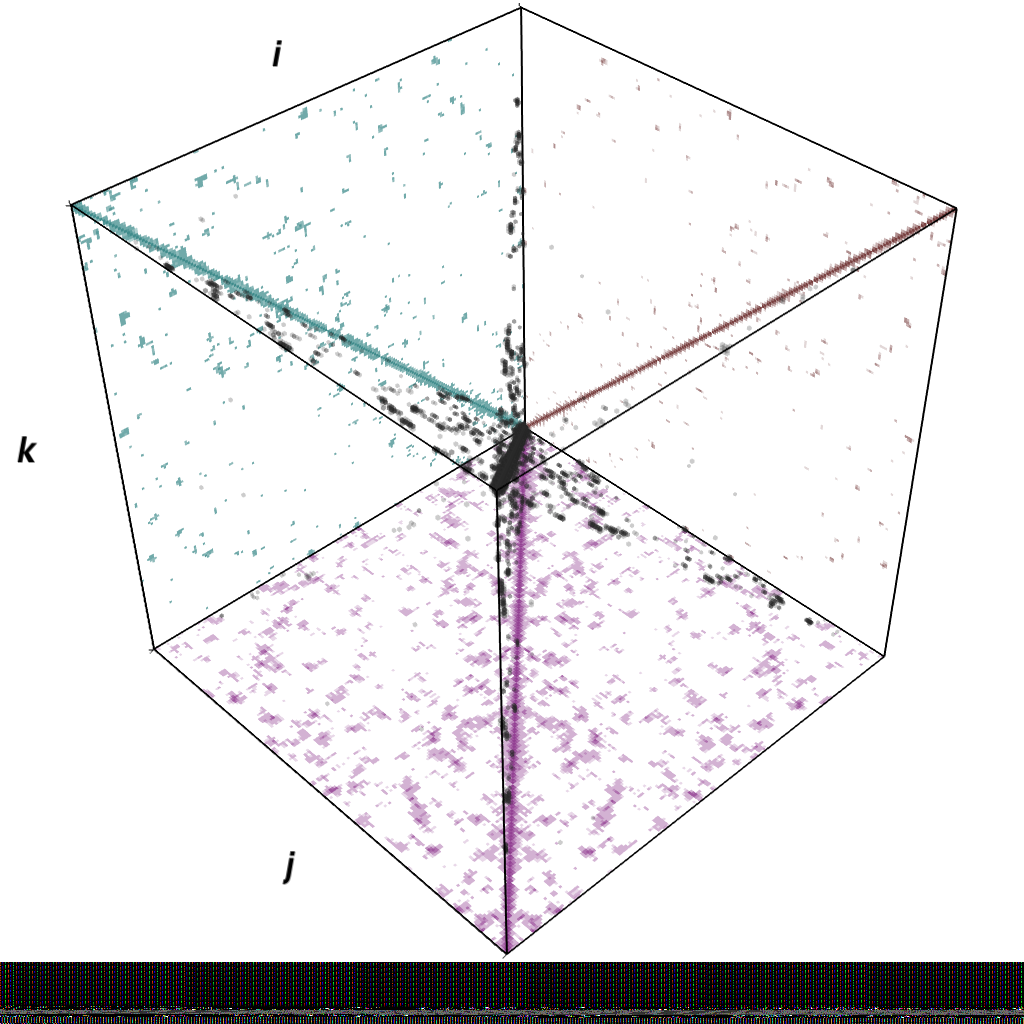}}
\fbox{ \includegraphics[width=3.8cm,keepaspectratio=true,
                        trim={0.cm 2.3cm 2.cm 1.cm},clip]
                        {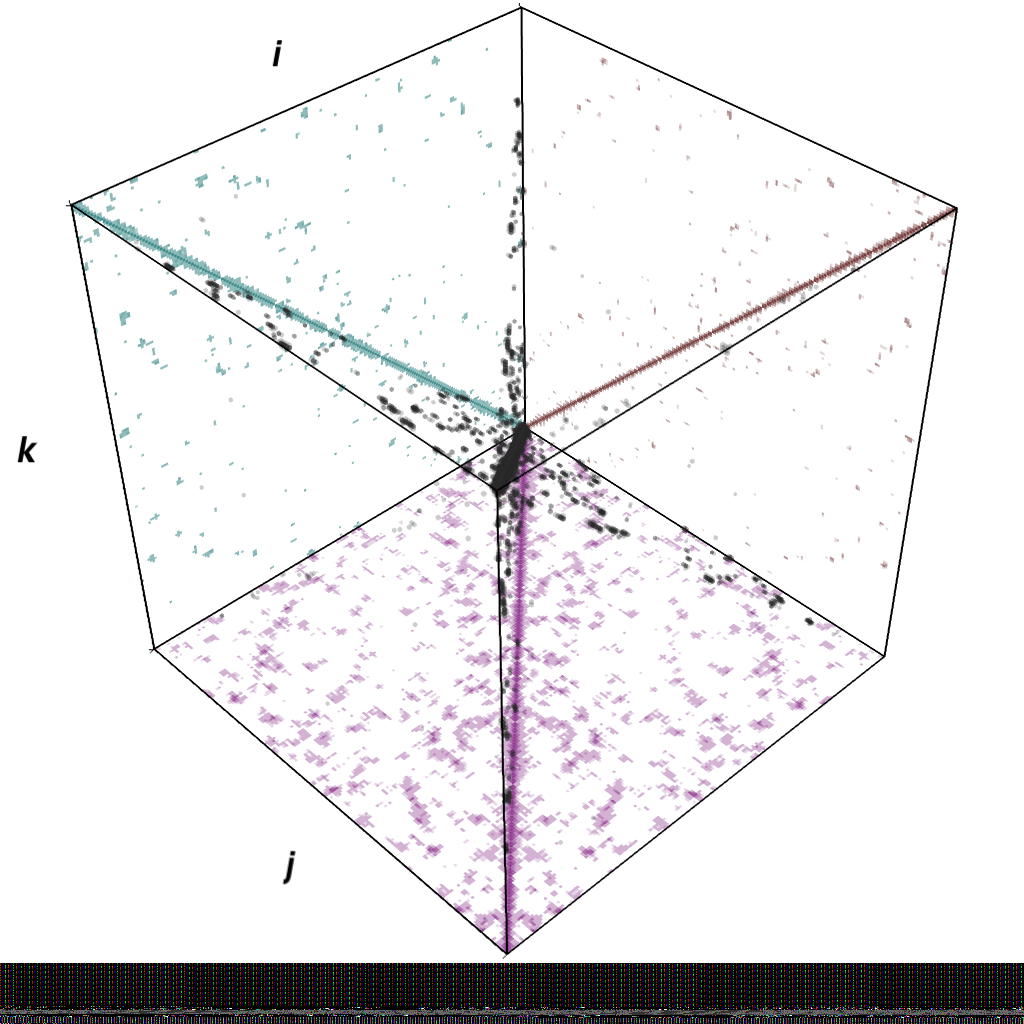}}
\caption{Product volumes in construction of the MAYEBOO preconditioner
$\left|\tau_0=.1,\mu_0=.1\, ; \,\scriptstyle{\mat{s}^{-1/2}} \right>$, with
dual instance square root iteration, for 6-311G** box of 100 periodic water molecules.
At top its  $\mat{y}_k=h_\alpha[ \mat{x}_{k-1} ] \ots \mat{y}_{k-1}$
for $k=0,4,\& 15$, while on the bottom we have $\mat{x}_k=  \mat{y}_{k}  \ot \mat{z}_{k}$ for $k=0,4, \& 15$.
Maroon is $\mat{a}$, purple is $\mat{b}$, green is $\mat{c}$,  and black is the volume ${\rm vol}_{a \ot b}$
in the product $\mat{c}=\mat{a} \ot \mat{b}$.}\label{Lensing2}
\end{figure}

In addition to compression through orderings that maximize these block-by-magnitude effects,
we demonstrate a new kind of locality in Figs.~\ref{Lensing1} and \ref{Lensing2},
which is, so far, uniquely exploited by the $n$-body approach to square root iteration.
This locality increases compressively towards convergence,
as contractive identity iterations develop.
We call this compression {\bf  \em lensing},  involving collapse of the culled volume about plane diagonals of the identity.
Lensing corresponds to strengthening Eq.(\ref{bound}), viz Eqs.~(\ref{boundY})-(\ref{boundZ}),
and to strongly contracting directional derivatives, viz Eqs.~(\ref{yorbit})-(\ref{zorbit}).
This is an important, mitigating effect for $\tt SpAMM$ computations in the $\mat{y}$ channel, encumbered by
the parameter $\tau_s \sim {\tt .01} \times \tau$.

\begin{figure}[h]
\fbox{ \includegraphics[width=6cm,keepaspectratio=true,
                        trim={0.cm 2.3cm 2.cm 1.cm},clip]
                        {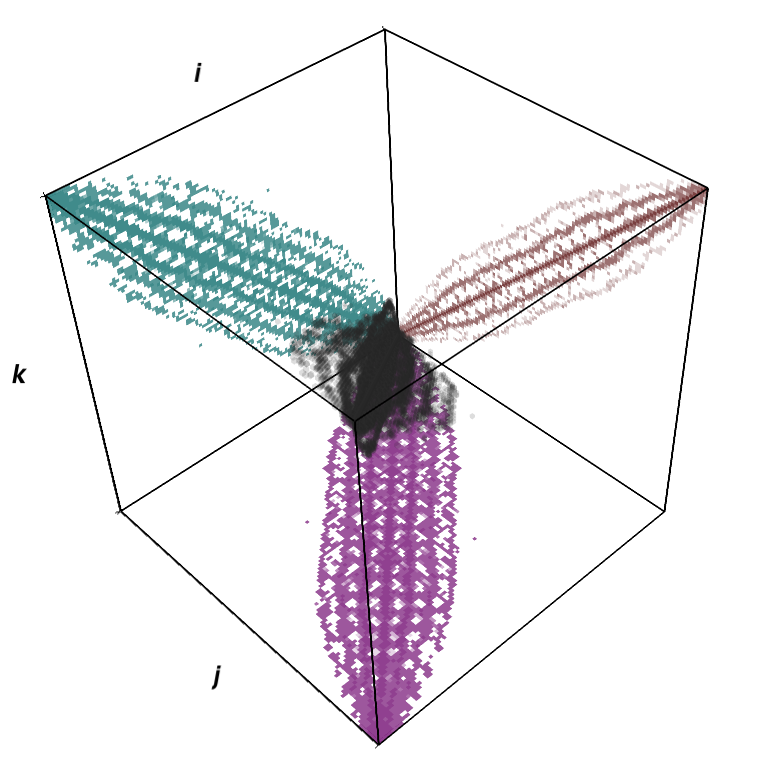}}
\fbox{ \includegraphics[width=6cm,keepaspectratio=true,
                        trim={0.cm 2.3cm 2.cm 1.cm},clip]
                        {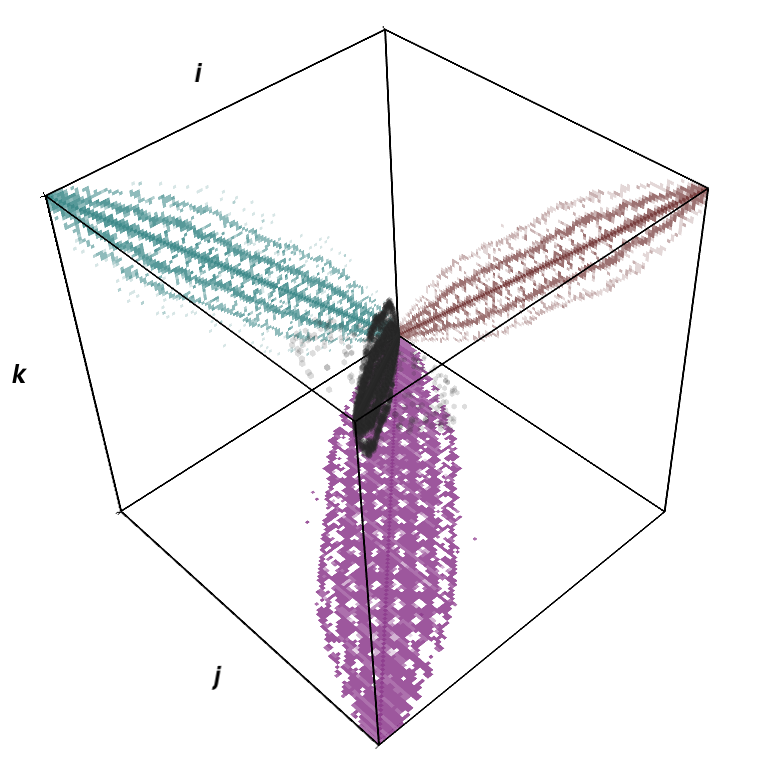}}

\fbox{ \includegraphics[width=6cm,keepaspectratio=true,
                        trim={0.cm 2.3cm 2.cm 1.cm},clip]
                        {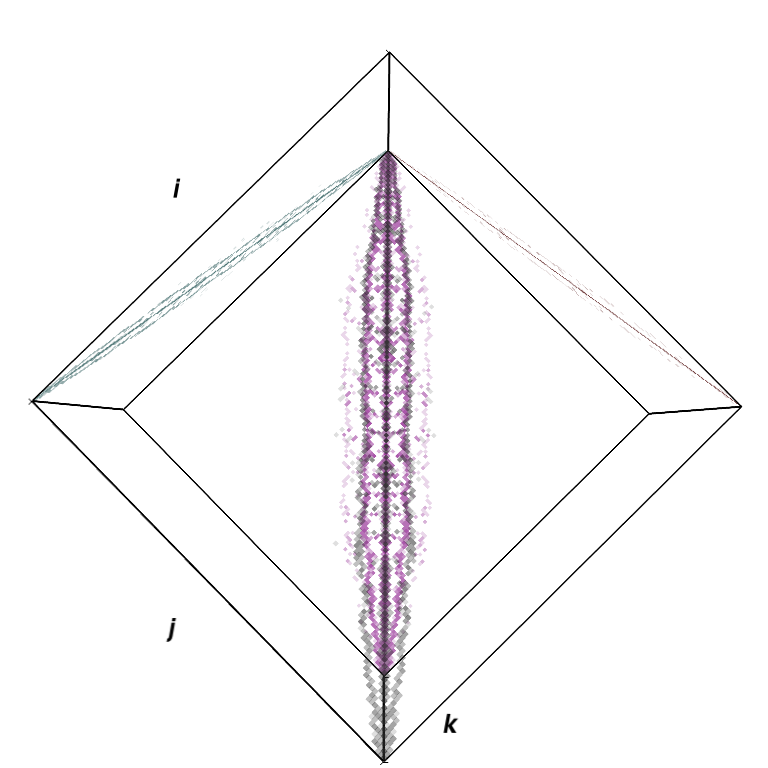}}
\fbox{ \includegraphics[width=6cm,keepaspectratio=true,
                        trim={0.cm 2.3cm 2.cm 1.cm},clip]
                        {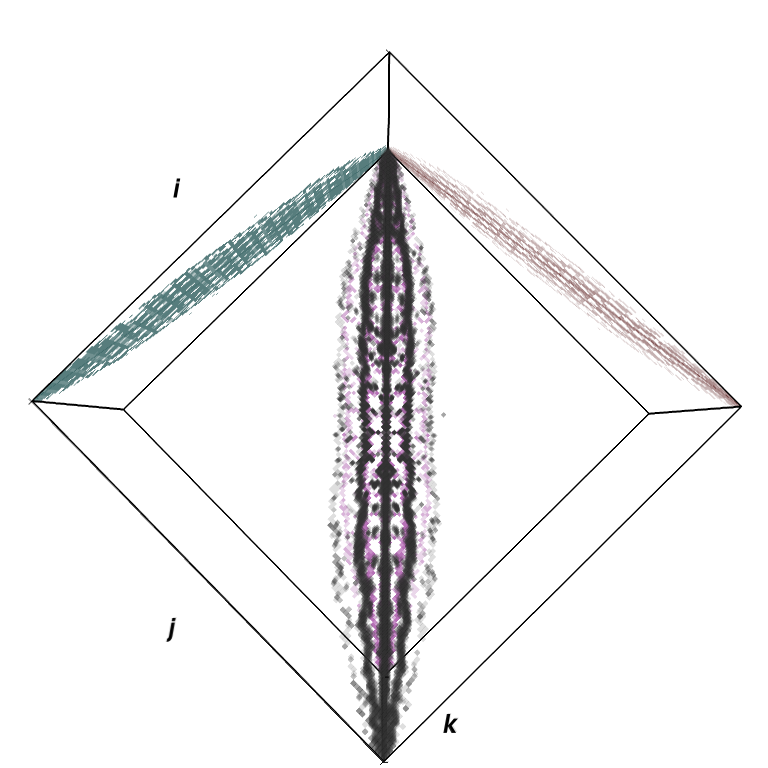}}
\caption{
Product volumes in construction of the unregularized preconditioner
$\left|\tau_0=.001,\mu_0=.0\, ; \,\scriptstyle{\mat{s}^{-1/2}} \right>$, with the
dual instance of square root iteration  and for the {\tt bcsstk14} structural matrix.
At top its  $\mat{y}_k=h_\alpha[ \mat{x}_{k-1} ] \ots \mat{y}_{k-1}$
for $k=0 \; \& \; 37$, while on the bottom we have $\mat{x}_k=  \mat{y}_{k}  \ot \mat{z}_{k}$ for $k=0\; \& \; 37$.
Maroon is $\mat{a}$, purple is $\mat{b}$, green is $\mat{c}$,  and black is the volume ${\rm vol}_{a \ot b}$
in the product $\mat{c}=\mat{a} \ot \mat{b}$.}\label{Lensing4}
\end{figure}

Graph reorderings that minimize the distance of matrix elements from the diagonal also lead to matrix locality (aforementioned).
In Fig.~\ref{Lensing4} we show convergence of an unregularized (MAYSS) preconditioner for this type of ordering and the ${\tt bcsstk14}$ \cite{BCSSTK14}
structural matrix of the Omni Coliseum in Atlanta, with $\kappa(\mat{s})=10^{10}$.
These results show remarkable gossamer sheeting and flattening along plane diagonals, in Fig.~\ref{Lensing4}, at top for development of $\mat{y}_k$,
as well as hollow accumulation of ${\rm vol}_{\mat{y}_k \ot \mat{z}_k}$ at bottom.
Interestingly, this example demonstrates well lensed volumes towards convergence, whilst the
equally ill-conditioned {and} lower dimensional  $\kappa (\mat{s})=10^{10}$ nanotube
demands a much tighter value of $\tau_s$ ($10^{-4}$ {vs.}~$10^{-9}$) and retains dense volumes through $128 \times$ the {unit cell}.

\begin{figure}[b]
\includegraphics[width=5.0in,keepaspectratio=true,trim={0.cm 0.cm 0.cm 0.cm},clip]
                 {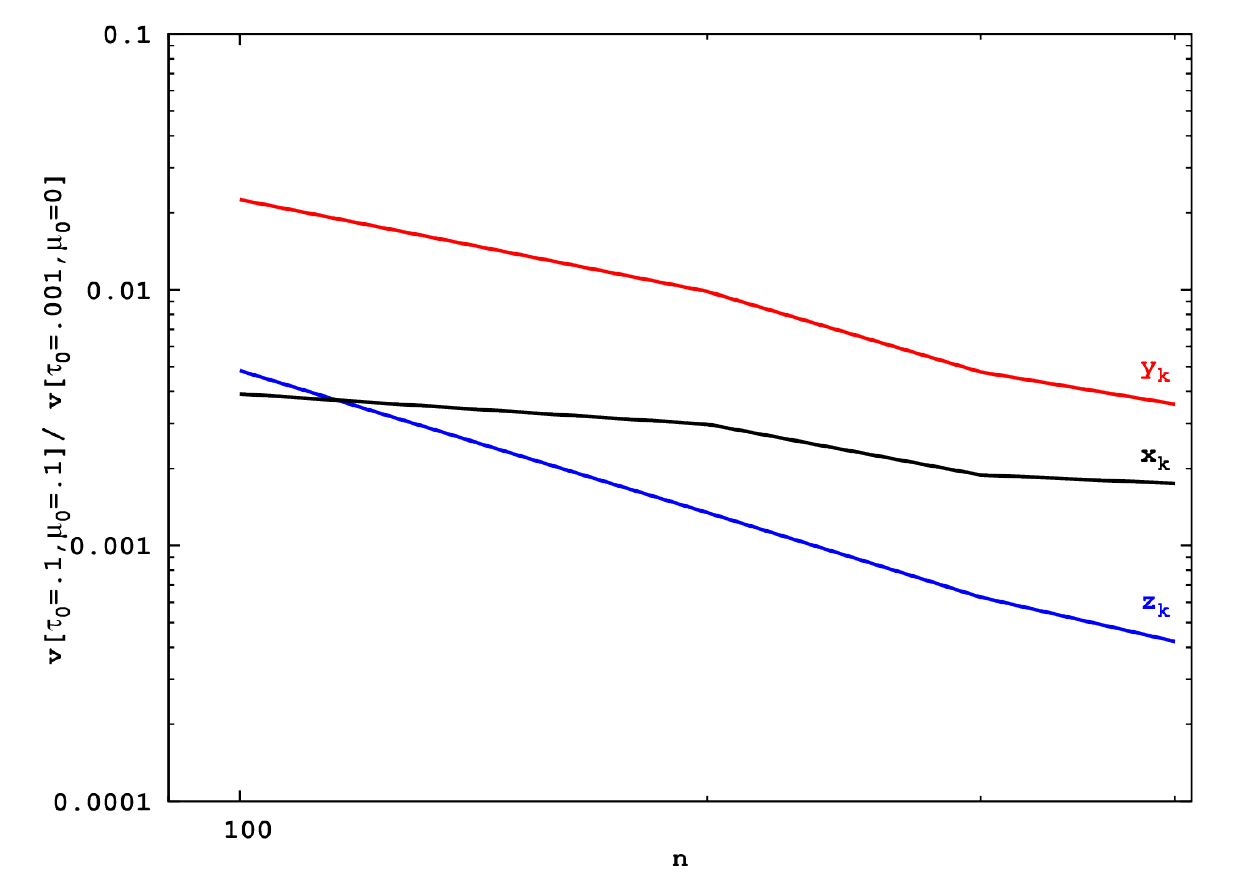}
\caption{
Complexity reduction in metric square root iteration for the periodic 6-311G** water sequence.
Shown is the ratio of lensed product volumes for the regularized MAYEBOO
approximation with respect to the unregularized (MAYSS) approximation.}\label{Complex1}
\end{figure}

\begin{figure}[b]
\includegraphics[width=5.0in,keepaspectratio=true,trim={0.cm 0.cm 0.cm 0.cm},clip]
                 {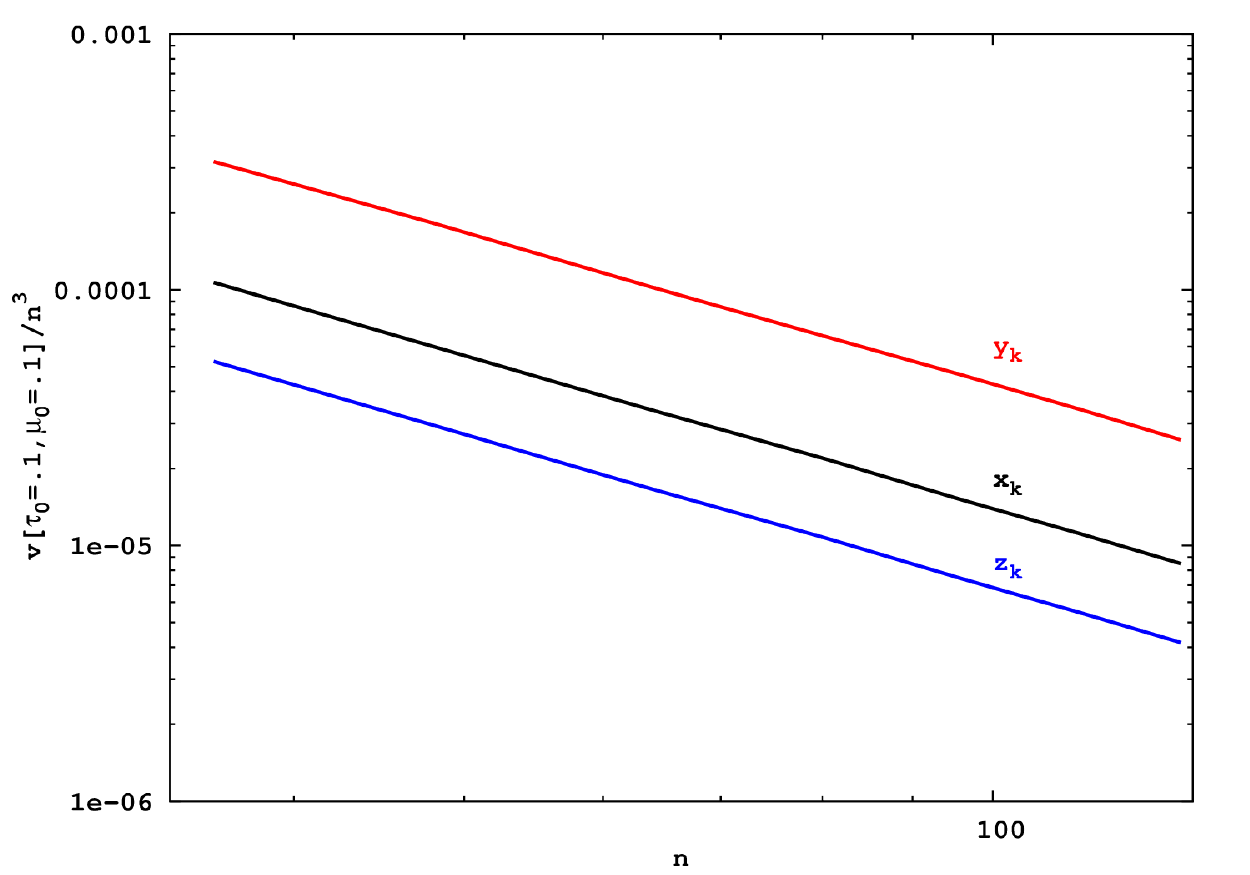}

\caption{
Complexity reduction in square root iteration for the $\kappa(\mat{s})=10^(10)$ sequence.
Shown is the ratio of lensed product volumes for the regularized MAYEBOO
approximation with respect to the unregularized MAYSS approximation, which we take to be $n^3$.}\label{Complex2}
\end{figure}

\subsection{Complexity reduction}

Finally, we show complexity reduction at convergence of the MAYEBOO approximation relative to the MAYSS approximation,
in Fig.~\ref{Complex1} for periodic water boxes, and in Fig.~\ref{Complex2} for the ill-conditioned nano-tube.
The two-orders difference between $\mat{y}_k$ and $\mat{z}_k$ volumes corresponds precisely to $\tau_s \sim \tau \times {\tt .01}$,
with $\mat{x}_k$ in between. Except for the slower trend in Fig.~(\ref{Complex1})'s $\mat{x}_k$ volume, we see the
potential for continued strong acceleration with increasing system size.

\section{Conclusions}\label{conclusions}

In this work, we developed the $\tt SpAMM$ $n$-body solver for square root iteration, along with
some algebra for the operator $\ot$, and showed how we could exploit different types of locality in the
sub-space metric of the product-tensor.
Our main contributions include a modified Cauchy-Schwarz criterion for the $\tt SpAMM$ occlusion-cull,
Eq.~(\ref{newspamm}), and proof that the corresponding relative error in the product is
bound by Eq.~(\ref{bound}).  We showed how block-by-magnitude orderings and locality
of the sub-space metric leads to reduced complexity of the $\tt SpAMM$ kernel, involving low-dimensional
sub-structures that bound the relative error,  distributed along plane-diagonals and along
their their intersection at the cube-diagonal.  Perhaps most significantly, we demonstrated a new kind of
compressive locality, lensing, that develops in the $\ot$ volume on contractive identity iteration,
together with tightening the $\tt SpAMM$ bound, {viz}~Eqs.~(\ref{boundY})-(\ref{boundZ}).

Additional contributions include development and implementation of a first order {Fr\'{e}chet} analyses for the single and dual
instances of the NS square root iteration, with focus on separating directional effects that are mostly {controlled} by the unperturbed
 reference algebra, from the magnitude of $\tt SpAMM$ occlusion errors and their accumulation.
We found that numerical sensitivity develops primarily in the $\mat{z}$ channel, according to Eq.~(\ref{zdispalcementbound}), due
to amplification of $\delta y$ by terms approaching condition of the full inverse;  we then looked at sensitivity to
this error, bifurcations, controled by $\tau_s$ (Figs.~\ref{flow_noscale_dual}-\ref{flow_scaled_stab}),
concluding that a most approximate, naive application of $\tt SpAMM$ to the ill-conditioned problem is generally insufficient
to achieve a fast solution.

Finally, we introduced scoping on both precision and regularization in product representation of the inverse factor,
and demonstrated the potential for orders of magnitude compression in the dual instance, Figs.~\ref{Complex1}-\ref{Complex2},
with the most extreme, ``by-one-order'' slice of the nested factor, providing a foothold for
this expansion at $\tau_0={\tt .1}$.   A next step is to demonstrate full bootstrapping of the inverse factor with reduced complexity,
{\em i.e.} via a compact, nested product of well lensed terms, a work in progress.

\doingrevtex{
\begin{acknowledgments}
This article was released under LA-UR-15-26304.  The Los Alamos National
Laboratory is operated by Los Alamos National Security, LLC for the NNSA of the
USDoE under Contract No.  DE-AC52- 06NA25396.
\end{acknowledgments}
}

\bibliographystyle{siam}

\bibliography{challacombe_haut_bock_2015_nbody_square_root_iteration}

\end{document}